\documentclass[aos,seceqn,dvips]{arximspdf}
\usepackage{mathbh}

\doi{10.1214/009053606000000614}
\volume{34}
\issue{4}
\pubyear{2006}
\firstpage{2026}
\lastpage{2068}

\makeatletter
\newproclaim{exa}{Example}
\newproclaim{definition}{Definition}
\newproclaim{ass}{Assumption}
\newtheorem{corollary}{Corollary}
\newproclaim{rem}{Remark}
\newtheorem{lemma}{Lemma}
\newtheorem{theorem}{Theorem}
\makeatother

\begin{document}
\begin{frontmatter}

\title{Efficient likelihood estimation in state~space~models\thanksref{T1}}
\runtitle{Likelihood estimation in SSM}

\begin{aug}
\author[A]{\fnms{Cheng-Der} \snm{Fuh}\corref{}\ead[label=e1]{stcheng@stat.sinica.edu.tw}}
\runauthor{C.-D. Fuh}
\thankstext{T1}{Supported in part by NSC 95-2118-M-001-008.}
\affiliation{National Central University and Academia Sinica}
\address[A]{Institute of Statistical Science\\
Academia Sinica\\
Taipei 11529\\
Taiwan\\
Republic of China\\
\printead{e1}} 
\end{aug}

\received{\smonth{3} \syear{2004}}
\revised{\smonth{9} \syear{2005}}

\begin{abstract}
Motivated by studying asymptotic properties of the maximum likelihood
estimator (MLE) in stochastic volatility (SV) models, in this paper
we investigate likelihood estimation in state space models. We first
prove, under some regularity conditions, there is a consistent sequence
of roots of the likelihood equation that is asymptotically normal with
the inverse of the Fisher information as its variance. With an extra
assumption that the likelihood equation has a unique root for each $n$,
then there is a consistent sequence of estimators of the unknown
parameters. If, in addition, the supremum of the log likelihood
function is integrable, the MLE exists and is strongly consistent.
Edgeworth expansion of the approximate solution of likelihood equation
is also established. Several examples, including Markov switching
models, ARMA models, (G)ARCH models and
stochastic volatility (SV) models, are given for illustration.
\end{abstract}

\pdfkeywords{62M09, 62F12, 62E25, Consistency,
efficiency, ARMA models, (G)ARCH models,
stochastic volatility models, asymptotic normality,
asymptotic expansion, Markov switching models,
maximum likelihood, incomplete data, iterated random functions,}

\end{frontmatter}

\section{Introduction}\label{s1}
Motivated by studying asymptotic properties of the maximum likelihood
estimator (MLE) in stochastic volatility (SV) models, in this paper
we investigate likelihood estimation in state space models. A state
space model is, loosely speaking, a sequence~$\{\xi_n\}_{n=0}^{\infty}$
of random variables obtained in the following way. First, a realization
of a Markov chain $\mathbf{X}=\{X_n, n \geq 0 \}$ is created. This
chain is sometimes called the regime and is not observed. Then,
conditional on $\mathbf{X}$, the $\xi$-variables are generated.
Usually the dependence of $\xi_n$ on $\mathbf{X}$ is more or less
local, as when $\xi_n = g(X_n,\xi_{n-1},\eta_n)$ for some function $g$
and random sequence~$\{\eta_n \}$, independent of $\mathbf{X}$. $\xi_n$
itself is generally not Markov and may, in fact, have a complicated
dependence structure. When the state space of $\{X_n, n \geq 0 \}$ is
finite, it is the so-called hidden Markov model or Markov switching
model.

The statistical modeling and computation for state space models have
attracted a great deal of attention recently because of their importance
in applications to speech recognition~\cite{RJ93}, signal processing
\cite{EAM95}, ion channels \cite{BR92}, molecular
biology~\cite{KBMSH94} and economics~\cite{B86,E82,T86}. The reader is
referred to~\cite{FY03,H94,K01} for a comprehensive summary. The main
focus of these efforts has been state space modeling and estimation,
algorithms for fitting these models and the implementation of
likelihood based methods.

The state space model here is defined in a general sense, in which the
observations are \textit{conditionally Markovian dependent}, and the
state space of the driving Markov chain need not be \textit{finite}
or \textit{compact}. When the state space is finite and
the observation is a deterministic function of the state space, Baum
and Petrie \cite{BP66} established the consistency and asymptotic
normality of the MLE. When the observed random variables are
conditionally independent, Leroux \cite{L92} proved strong consistency
of the MLE, while Bickel, Ritov and Ryd{\'e}n \cite{BRR98} established
asymptotic normality of the MLE under mild conditions. Jensen and
Petersen \cite{JP99}, Douc and Matias \cite{DM01} and Douc, Moulines
and Ryd{\'e}n \cite{DMR04} studied asymptotic properties of the MLE for
general ``pseudo-compact'' state space models. By extending the
inference problem to time series analysis where the state space is
finite and the observed random variables are conditionally Markovian
dependent, Goldfeld and Quandt \cite{GQ73} and Hamilton \cite{H89}
considered the implementation of the maximum likelihood estimator in
switching autoregressions with Markov regimes. Francq and Roussignol
\cite{FR98} studied the consistency of the MLE, while Fuh \cite{F04a}
established the Bahadur efficiency of the MLE in Markov switching
models. We now give two examples of state space models.

\begin{exa}[{[$\mbox{GARCH}(p,q)$ model]}]\label{e1}
For given $p \geq 1$ and $q\geq 0$, let
\begin{equation}\label{6.8}
Y_n = \sigma_n \varepsilon_n\quad\mbox{and}\quad \sigma_n^2 = \delta +
\sum_{i=1}^p \alpha_i \sigma_{n-i}^2 + \sum_{j=1}^q \beta_j Y_{n-j}^2,
\end{equation}
where $\delta > 0,  \alpha_i \geq 0$  and $\beta_j \geq 0$ are
constants, $\varepsilon_n$ is a sequence of independent and identically
distributed (i.i.d.) random variables, and $\varepsilon_n$ is
independent of $\{Y_{n-k}, k \geq 1\}$ for all $n$. This is the
celebrated $\mbox{GARCH}(p,q)$ model proposed by Bollerslev
\cite{B86}. When $q=0$ or $\beta_j=0$, for $j=1,\ldots,q$, this is the
$\mbox{ARCH}(p)$ model first considered by Engle \cite{E82}. The
reader is referred to~\cite{BEN94}~and~\cite{FY03} for a comprehensive
summary.

For convenience of notation, we assume that $p,q \geq 2$, and by adding
some $\alpha_i$ or $\beta_j$ equal to zero if necessary.
Denote
$\eta_n = \sigma_n^{-1} Y_n$, $\tau_n = (\alpha_1 + \beta_1 \eta_n^2,
\alpha_2,\ldots, \alpha_{p-1}) \in \mathbf{R}^{p-1}$, $\zeta_n =
(\eta_n^2,0,\ldots,0) \in \mathbf{R}^{p-1}$, $\beta =
(\beta_2,\ldots,\beta_{q-1}) \in \mathbf{R}^{q-2}$, and let $I_{p-1}$
and $I_{q-2}$ be identity matrices. Let $A_n$ be a $(p+q-1) \times
(p+q-1)$ matrix written in block form as
\begin{equation}\label{6.10}
A_{n} =  \left[\matrix{
 \tau_n & \alpha_p & \beta & \beta_q \cr\noalign{}
 I_{p-1} & 0 & 0 & 0 \cr\noalign{}
 \zeta_n & 0 & 0 & 0 \cr\noalign{}
 0 & 0 & I_{q-2} & 0}\right].
\end{equation}
Note that $\{A_n, n \geq 0 \}$ are i.i.d. random matrices.

Let $Z=(\delta,0,\ldots,0)' \in \mathbf{R}^{p+q-1}$ and
$X_n=(\sigma_{n+1}^2,\ldots,\sigma_{n-p+2}^2,
Y_n^2,\ldots,\break Y_{n-q+2}^2)'$, where ``$\,'$'' denotes transpose. Following  the
idea of Bougerol and\break \mbox{Picard~\cite{10}}, we have the following state
space representation of the\break $\mbox{GARCH}(p,q)$ model: $X_n$ is a
Markov chain governed by
\begin{equation}\label{6.11}
X_{n+1} = A_{n+1} X_n + Z,
\end{equation}
and $\xi_n := g(X_n)=(Y_n^2,\ldots,Y_{n-q+2}^2)'$, the observed random
quantity, is a noninvertible function of $X_n$.
\end{exa}

\begin{exa}[(Stochastic volatility models)]\label{e2}
Let
\begin{equation}\label{6.12}
Y_n = \sigma_n \varepsilon_n,
\end{equation}
where $\log \sigma_n^2$ follows an $\mbox{AR}(1)$ process and $\varepsilon_n$
is a sequence of i.i.d. random variables with standard normal
probability density function. This is the discrete time stochastic
volatility model proposed by Taylor \cite{T86}. The reader is referred
\mbox{to \cite{CHR96,S96,T94}} for a comprehensive summary. Note that
Genon-Catalot, Jeantheau and Lar\'{e}do \cite{GJL00} studied the
ergodicity and mixing properties of stochastic volatility models from
the hidden Markov model point of view.

Write $X_n := \log \sigma_n^2$ and $Y_n = \sigma \varepsilon_n
\exp(X_n/2),$ where $\sigma$ is a scale parameter. Squaring the
observations in the above equation and taking logarithms gives $\log
Y_n^2 = \log \sigma^2 +  X_n + \log \varepsilon_n^2.$ Alternatively, we
have
\begin{equation}\label{6.15}
\log Y_n^2 = \omega +  X_n + \zeta_n,
\end{equation}
where $\omega = \log \sigma^2 + E \log \varepsilon_n^2,$ so that the
disturbance $\zeta_n$ has mean zero by construction. The scale
parameter $\sigma$ also removes the need for a constant term in the
stationary first-order autoregressive process
\begin{equation}\label{6.16}
X_{n} = \alpha X_{n-1} + \eta_n,\qquad |\alpha| < 1,
\end{equation}
where $\eta_n$ is a sequence of i.i.d. random variables
distributed as $N(0,\sigma_{\eta}^2)$. Moreover, we assume that $\zeta_n$ and $\eta_n$
are independent. Note that in (\ref{6.15}) and (\ref{6.16}) the
observed random quantity is $\xi_n:=\log Y_n^2$. $\{X_n, n \geq 0\}$
and forms a Markov chain with transition probability
\begin{equation}\label{6.17}
p(x_{k-1},x_k) = (2 \pi \sigma_\eta^2)^{-1/2} \exp  \biggl\{ -
\frac{1}{2} \frac{( x_k - \alpha x_{k-1})^2}{\sigma_\eta^2} \biggr\}
\end{equation}
and stationary distribution $\pi \sim N(0,\sigma_\eta^2/(1-\alpha))$.

For given observations $\mathbf{ y}=(\log y_1^2,\ldots,\log y_n^2)$
from the state space model (\ref{6.15}) and (\ref{6.16}), the
likelihood function of the parameter $\theta=(\alpha,\sigma_\eta^2)$ is
\begin{eqnarray}\label{6.18}
l(\mathbf{y};\theta) &=& \int_{x_0 \in {\mathcal  X}} \cdots \int_{x_n
\in {\mathcal  X}} \pi(x_0) \prod_{k=1}^n p(x_{k-1},x_k) \nonumber
\\[-8pt]
\\[-8pt]
\nonumber &&\hspace*{39mm}{}\times   f_{\zeta}(\log y_k^2 - \omega -
x_{k})\,dx_n \cdots dx_0,
\end{eqnarray}
where $f_{\zeta}(\cdot)$ is the probability density function of
$\zeta_1$.

A major difficulty in analyzing the likelihood function in state space
models is that it can be expressed only in integral form; see equation
(\ref{6.18}), for instance. In this paper we provide a device which
represents the integral likelihood function as the $L_1$-norm of a
Markovian iterated random functions system. This new representation
enables us to apply results of the strong law of large numbers, central
limit theorem and Edgeworth expansion for the distributions of Markov
random walks, and to verify strong consistency of the MLE and first-order
efficiency and Edgeworth expansion on the solution of the likelihood
equation. Note that third-order efficiency follows from Edgeworth
expansion by a standard argument (cf.~\cite{G94}). Another essential
point worth being mentioned is that we introduce \textit{a weight
function} in a suitable way [see (\ref{4.1})--(\ref{4.3}),
Assumptions \hyperref[assK2]{K2}, \hyperref[assK3]{K3} and Definition~\ref{D2} in
Section~\textup{\ref{s4}}, and C1 in Section~\textup{\ref{s5}}] to
relax the condition of a compact state space for the underlying Markov
chain, and to cover several interesting examples.

The remainder of this paper is organized as follows. In
Section~\textup{\ref{s2}} we define the state space model as a general
state Markov chain in a Markovian random environment, and represent the
likelihood function as the $L_1$-norm of a Markovian iterated random
functions system. In Section~\textup{\ref{s3}} we give a brief summary
of a Markovian iterated random functions system, and provide an  ergodic
theorem and the strong law of large numbers. The multivariate central
limit theorem and Edgeworth expansion for a Markovian iterated random
functions system are given in Section~\textup{\ref{s4}}.
Section~\textup{\ref{s5}} contains our main results, where we
consider efficient likelihood estimation in state space models, and
state the main results. First, we compute Fisher information and
prove the existence of an efficient estimator in a ``Cram{\'e}r
fashion.'' Second, we characterize Kullback--Leibler information,
and prove strong consistency of the MLE. Last, we establish Edgeworth
expansion of the approximate solution of the likelihood equation. In
Section~\textup{\ref{s6}} we consider a few examples, including Markov
switching models, \textup{ARMA} models, (G)\mbox{ARCH} models and
SV models, which are commonly used in financial economics. The proofs
of the lemmas in Section~\textup{\ref{s5}} are given in
Section~\textup{\ref{s7}}. Other technical proofs are deferred to the
\hyperref[app]{Appendix}.
\end{exa}

\section{State space models}\label{s2}
A state space model is defined as a
parameterized Markov chain in a Markovian random environment with the
underlying environmental Markov chain viewed as missing data.
Specifically, let $\mathbf{X}= \{X_n, n \geq 0 \}$ be a Markov chain on
a general state space ${\mathcal  X}$, with transition probability
kernel $P^{\theta}(x,\cdot)= P^{\theta}\{X_1 \in \cdot|X_0=x\}$ and
stationary probability $\pi_{\theta}(\cdot)$, where $\theta \in \Theta
\subseteq \mathbf{R}^q $ denotes the unknown parameter. Suppose that a
random sequence $\{\xi_n\}_{n=0}^{\infty},$ taking values in
$\mathbf{R}^d$, is adjoined to the chain such that $\{(X_n,\xi_n), n
\geq 0\}$ is a Markov chain on ${\mathcal  X} \times \mathbf{R}^d$
satisfying $P^{\theta} \{ X_1 \in A | X_0=x,\xi_0=s \} = P^{\theta} \{
X_1 \in A | X_0=x \}$ for $A \in {\mathcal  B}({\mathcal  X})$, the
$\sigma$-algebra of ${\mathcal X}$. And conditioning on the full
$\mathbf{X}$ sequence, $\xi_n$~is a Markov chain with probability
\begin{eqnarray}\label{2.1}
&&  P^{\theta} \{\xi_{n+1} \in  B | X_0,X_1,\ldots;\xi_0,\xi_1,\ldots, \xi_n \}
\nonumber
\\[-8pt]
\\[-8pt]
\nonumber
&&\qquad = P^{\theta} \{\xi_{n+1} \in  B | X_{n+1};\xi_n \}\qquad\mbox{a.s.}
\end{eqnarray}
for each $n$ and $B \in {\mathcal  B}(\mathbf{R}^d),$ the Borel
$\sigma$-algebra on $\mathbf{R}^d$. Note that in (\ref{2.1}) the
conditional probability of $\xi_{n+1}$ depends on $X_{n+1}$ and $\xi_n$
only. Furthermore, we assume the existence of a transition probability
density $p_{\theta}(x,y)$ for the Markov chain $\{X_n, n \geq 0\}$ with
respect to a $\sigma$-finite measure $m$ on ${\mathcal  X}$ such that
\begin{eqnarray}\label{2.2}
&&   P^{\theta} \{X_1 \in A, \xi_{1} \in  B | X_0=x, \xi_0 =s_0   \}
\nonumber \\[-8pt]
\\[-8pt]
\nonumber &&\qquad = \int_{y \in A} \int_{s \in B} p_{\theta}(x,y)
f(s;\theta|y, s_0) Q(ds) m(dy),
\end{eqnarray}
where $f(\xi_k;\theta|X_k,\xi_{k-1})$ is the conditional probability
density of $\xi_k$ given $\xi_{k-1}$ and $X_k$, with respect to a
$\sigma$-finite measure $Q$ on $\mathbf{R}^d$. We also assume that the
Markov chain $\{(X_n,\xi_n), n \geq 0\}$ has a stationary probability
with probability density function $\pi(x)f(\cdot;\theta|x)$ with
respect to $m \times Q$. In this paper we consider
$\theta=(\theta_1,\ldots,\theta_q) \in \Theta \subseteq \mathbf{R}^q$
as the unknown parameter, and the true parameter value is denoted by
$\theta_0$. We will use $\pi(x)$ for $\pi_{\theta}(x)$, $p(x,y)$ for
$p_{\theta}(x,y)$, $f(\xi_0|X_0)$ for $f(\xi_0;\theta|X_0)$, and
$f(\xi_k|X_k, \xi_{k-1})$ for $f(\xi_k;\theta|X_k,\xi_{k-1})$,
here and in the sequel, depending on our convenience. Now
we give a formal definition as follows.\looseness=1

\begin{definition}\label{def1}
$\{\xi_n,n \geq 0\}$ is called a state space model if there is a Markov
chain $\{X_n,n \geq 0\}$ such that the process $\{(X_n,\xi_n),n \geq
0\}$ satisfies~(\ref{2.1}).
\end{definition}

Note that this setting includes several interesting examples of
Markov-switching Gaussian autoregression of Hamilton \cite{H89},
(G)\mbox{ARCH} models of Engle \cite{E82} and Bollerslev
\cite{B86}, and SV models of Clark \cite{C73} and Taylor \cite{T86}.
When the state space ${\mathcal  X}$ is finite or compact, this reduces
to the hidden Markov model considered by Francq and Roussignol
\cite{FR98}, Fuh \cite{F03,F04a,F04c} and Douc, Moulines and Ryd{\'e}n
\cite{DMR04}. Denote $S_n=\sum_{t=1}^n \xi_t$. When $\xi_n$ are
\textit{conditionally independent} given $\mathbf{X}$, the Markov chain
$\{(X_n,S_n),n \geq 0\}$ is called a \textit{Markov additive process}
and $S_n$ is called a \textit{Markov random walk}. Furthermore, if the
state space ${\mathcal X}$ is finite, $\{\xi_n,n \geq 0\}$ is the
hidden Markov model studied by Leroux \cite{L92}, Bickel and Ritov
\cite{BR96} and Bickel, Ritov and Ryd{\'e}n \cite{BRR98}. When the
state space ${\mathcal X}$ is ``pseudo-compact'' and $\xi_n$ are
conditionally independent given $\mathbf{X}$, $\{\xi_n,n \geq 0\}$ is
the state space model considered
in~\cite{JP99}~and~\cite{DM01}.

For given observations $s_0,s_1,\ldots,s_n$ from a state space model
$\{\xi_n,n \geq 0\}$, the likelihood function is
\begin{eqnarray}\label{2.3}
\nonumber && p_n(s_0,s_1,\ldots,s_n;\theta)
\\
&&\qquad = \int_{x_0 \in {\mathcal  X}} \cdots  \int_{x_n \in {\mathcal
X}} \pi_\theta(x_0) f(s_0;\theta|x_0)
\nonumber
\\[-8pt]
\\[-8pt]
\nonumber
&&\hspace*{35mm} {}\times  \prod_{j=1}^n p_\theta(x_{j-1},x_j)
\\
\nonumber &&\hspace*{44.2mm}
{}\times f(s_j;\theta| x_j, s_{j-1})m(dx_n)\cdots m(dx_0).
\end{eqnarray}
Recall that $\pi_\theta(x_0) f(s_0;\theta|x_0)$ is the stationary
probability density with respect to $m \times Q$ of the Markov chain
$\{(X_n,\xi_n),n \geq 0\}$.

To represent the likelihood $p_n(\xi_0,\xi_1,\ldots,\xi_n;\theta)$ as
the $L_1$-norm of a Markovian iterated random functions system, let
\begin{equation}\label{2.4}
\hspace*{12pt}
\mathbf{M} = \biggl\{h|h\dvtx {\mathcal  X} \to \mathbf{R}^+  \mbox{ is
}m\mbox{-measurable and } \int_{x \in {\mathcal  X}} h(x)m(dx) < \infty
\biggr\}.
\end{equation}
For each $j=1,\ldots,n$, define the random functions
$\mathbf{P}_\theta(\xi_0)$ and $\mathbf{P}_\theta(\xi_j)$ on
$({\mathcal X} \times \mathbf{R}^d) \times \mathbf{M}$ as
\begin{eqnarray}
  \mathbf{P}_\theta(\xi_0)h(x) &=&  \int_{x \in {\mathcal  X}} f(\xi_0;\theta|x) h(x) m(dx),\qquad\mbox{a constant,}\label{2.5}
\\
  \mathbf{P}_\theta(\xi_j)h(x) &=&  \int_{y \in {\mathcal  X}}
p_{\theta}(x,y) f(\xi_j;\theta|y, \xi_{j-1}) h(y) m(dy). \label{2.6}
\end{eqnarray}
Define the composition of two random functions as
\begin{eqnarray}\label{2.7}
\nonumber &&
\mathbf{P}_\theta(\xi_{j+1})\circ\mathbf{P}_\theta(\xi_{j})h(x)
\\
&&\qquad = \int_{z \in {\mathcal  X}} p_\theta(x,z) f(\xi_j;\theta|z,
\xi_{j-1})
\\
\nonumber &&\hspace*{19mm}  {}\times \biggl(\int_{y \in {\mathcal  X}} p_\theta(z,y)
f(\xi_{j+1};\theta|y,\xi_{j}) h(y) m(dy) \biggr) m(dz).
\end{eqnarray}
For $h \in \mathbf{M}$, denote $\|h\|:= \int_{x \in {\mathcal  X}} h(x)
m(dx)$ as the $L^1$-norm on $\mathbf{M}$ with respect to~$m$. Then the
likelihood $p_n(\xi_0,\xi_1,\ldots,\xi_n;\theta)$ can be represented as
\begin{eqnarray}\label{2.8}
\nonumber && p_n(\xi_0,\xi_1,\ldots,\xi_n;\theta)
\\
\nonumber &&\qquad = \int_{x_0 \in {\mathcal  X}} \cdots  \int_{x_n \in {\mathcal
X}} \pi_\theta(x_0) f(\xi_0;\theta|x_0)
\\
&&\hspace*{35.2mm} {}\times \prod_{j=1}^n p_\theta(x_{j-1},x_j)
\\
\nonumber &&\hspace*{44mm} {}\times f(\xi_j;\theta|x_j,\xi_{j-1})m(dx_n)\cdots m(dx_0)
\\
\nonumber &&\qquad =  \|\mathbf{P}_\theta(\xi_n) \circ \cdots \circ
\mathbf{P}_\theta(\xi_1) \circ \mathbf{P}_\theta(\xi_0) \pi_{\theta}\|.
\end{eqnarray}

Note that, for $j=1,\ldots,n,$ the integrand
$p_\theta(x,y)f(\xi_j;\theta|y,\xi_{j-1})$ of $\mathbf{P}_\theta(
\xi_j)$ in (\ref{2.6}) and (\ref{2.8}) represents $X_{j-1} =x$ and
$X_j\in dy$, and $\xi_j$ is a Markov chain with transition probability
density $f(\xi_j;\theta|y,\xi_{j-1})$ for given $\mathbf{X}$. By
definition (\ref{2.1}), $\{(X_n,\xi_n), n \geq 0 \}$ is a Markov chain,
and this implies that $\mathbf{P}_\theta(\xi_j)$ is a sequence of
Markovian iterated random functions systems (see
Section~\textup{\ref{s5}} for a formal definition). Therefore, by
representation (\ref{2.8}),
$p_n(\xi_0,\xi_1,\ldots,\xi_n;\theta)$ is the $L_1$-norm of a Markovian
iterated random functions system.

\section{Ergodic theorems for a Markovian iterated random functions system}\label{s3}
To analyze the asymptotic properties of efficient likelihood estimators
in state space models, in this section we study the ergodic theorem and the
strong law of large numbers for a Markovian iterated random functions
system. The Markovian iterated random functions system is a
generalization of an iterated random functions system, in which the
random functions are driven by a Markov chain. For a general account of
an iterated random functions system, the reader is referred
to~\cite{DF99} for a recent survey.

For simplicity in our notation, let $\{Y_n,n\geq 0\}$ [instead of
$\{(X_n,\xi_n), n \geq 0\}$ in Section~\textup{\ref{s2}}] be a Markov
chain on a general state space ${\mathcal  Y}$ with $\sigma$-algebra
${\mathcal A}$, which irreducible with respect to a maximal
irreducibility measure on $({\mathcal  Y},{\mathcal  A})$ and is
aperiodic. The transition kernel is denoted by $P(y,A)$. Let
$(\mathbf{M}, d)$ be a complete separable metric space with Borel
$\sigma$-algebra ${\mathcal B}(\mathbf{M})$. Denote by $M_0$ a random
variable which is independent of $\{Y_{n}, n \geq 0 \}$. A sequence of
the form
\begin{equation}\label{3.1}
 M_{n} = F(Y_{n},M_{n-1}),\qquad n\ge 1,
\end{equation}
taking values in $(\mathbf{M}, d)$ is called a \textit{Markovian
iterated random functions system} (MIRFS) of Lipschitz functions
providing the following:
\begin{longlist}[(2)]
\item[(1)] $\{Y_{n}, n \geq 0 \}$ is a Markov chain taking values in a
second countable measurable space $({\mathcal  Y},{\mathcal  A})$, with
transition probability kernel $P(\cdot,\cdot)$ and stationary
probability $\pi$, and $M_{0}$ is a random element on a probability
space $(\Omega,{\mathcal  F},P)$, which is independent of $\{Y_n,n \geq
0\}$;

\item[(2)] $F\dvtx ({\mathcal  Y} \times \mathbf{M},{\mathcal  A}
\otimes {\mathcal B}(\mathbf{M})) \to (\mathbf{M},{\mathcal
B}(\mathbf{M}))$ is
jointly measurable and Lipschitz continuous in the second argument.
\end{longlist}
Clearly, $\{(Y_n, M_n),n \geq 0\}$ constitutes a Markov chain with
state space ${\mathcal  Y} \times \mathbf{M}$ and transition
probability kernel $\mathbf{P}$, given by
\begin{equation}\label{3.2}
\mathbf{P}\bigl((y,u),A \times B\bigr) := \int_{z \in
A}I_B\bigl(F(z,u)\bigr) P(y,dz)
\end{equation}
for all $y \in {\mathcal  Y}, u \in \mathbf{M}, A \in {\mathcal  A}$
and $B \in {\mathcal  B}(\mathbf{M})$, where $I$ denotes the indicator
function. The $n$-step transition kernel is denoted $\mathbf{P}^{n}$.
For $(y,u) \in {\mathcal  Y} \times \mathbf{M}$, let $\mathbf{P}_{yu}$
be the probability measure on the underlying measurable space under
which $Y_0= y, M_{0}=u$ a.s. The associated expectation is denoted
$\mathbf{E}_{yu}$, as usual. For an arbitrary distribution $\nu$ on
${\mathcal Y} \times \mathbf{M}$, we put $\mathbf{P}_{\nu}(\cdot) :=
\int \mathbf{P}_{yu}(\cdot) \nu(dy \times du)$ with associated
expectation $\mathbf{E}_{\nu}$. We use $\mathbf{P}$ and $\mathbf{E}$
for probabilities and expectations, respectively, that do not depend on
the initial distribution.

Let $\mathbf{M}_{0}$ be a dense subset of $\mathbf{M}$ and ${\mathcal
M} (\mathbf{M_0},\mathbf{M})$ the space of all mappings $h\dvtx
\mathbf{M_0}\to\mathbf{M}$ endowed with the product topology and product
$\sigma$-algebra. Then the space ${\mathcal
L}_{\mathrm{Lip}}(\mathbf{M},\mathbf{M})$ of all Lipschitz continuous
mappings $h\dvtx {\mathbf M}\to{\mathbf M}$ properly embedded forms a
Borel subset of ${\mathcal  M}({\mathbf M_0},{\mathbf M})$, and the
mappings
\begin{eqnarray*}
{\mathcal  L}_{\mathrm{Lip}}({\mathbf M},{\mathbf M}) \times {\mathbf
M}  &\ni&  (h,u) \mapsto h(u) \in {\mathbf M},
\\
{\mathcal  L}_{\mathrm{Lip}}({\mathbf M},{\mathbf M})  &\ni& h \mapsto
l(h) := \sup_{u\ne v}{d(h(u),h(v))\over d(u,v)}
\end{eqnarray*}
are Borel; see Lemma 5.1 in~\cite{DF99} for details. Hence,
\begin{equation}\label{3.3}
L_{n} := l\bigl(F(Y_{n},\cdot)\bigr),\qquad  n\ge 0,
\end{equation}
are also measurable and form a sequence of Markovian dependent random
variables.

An important point to characterize the limit in the ergodic theorem will be
the right use of the idea of duality. For this purpose, we introduce a
time-reversed (or dual) Markov chain $\{\tilde{Y}_n, n \geq 0\}$ of
$\{Y_n, n \geq 0\}$ as follows. Assume that there exists a
\mbox{$\sigma$-}finite measure $m$ on $({\mathcal  Y},{\mathcal  A})$ such
that the probability measure $P$ on $({\mathcal  Y},{\mathcal  A})$
defined by $P(A)= P(Y_1 \in A|Y_0=y)$ is absolutely continuous with
respect to $m$, so that $P(A)= \int_A p(y,z)m(dz)$ for all $A\in
{\mathcal  A}$, where $p(y,\cdot)= dP/dm$. The Markov chain
$\{Y_n,n\geq 0\}$ is assumed to have an invariant probability measure
$\pi$ which has a positive probability density function $\pi$ (without
any confusion, we still use the same notation) with respect to $m$. We
shall use $\sim$ to refer to the time-reversed (or dual) process
$\{\tilde{Y}_n,n\geq 0\}$ with transition probability density
\begin{equation}\label{3.4}
 \tilde{p}(z,y)=p(y,z)\pi(y)/\pi(z).
\end{equation}
Denote $\tilde{P}$ as the corresponding probability. It is easy to see
that both $Y_n$ and $\tilde{Y}_n$ have the same stationary distribution
$\pi$. In this section we will assume that the initial distribution of
$Y_0$ is the stationary distribution~$\pi$.

In the following, we write $F_{n}(u)$ for $F(Y_{n},u)$. For all $1\le
k\le n$, let $F_{k:n} := F_{k}\circ \cdots \circ F_{n}$, $F_{n:k}:=
F_{n}\circ \cdots \circ F_{k}$, where $\circ$ denotes the composition
of functions. Denote $F_{n:n-1}$ as the identity on ${\mathbf M}$,
Hence
\begin{equation}\label{3.5}
 M_{n} = F_{n}(M_{n-1}) = F_{n:1}(M_{0})
\end{equation}
for all $n\ge 0$. Closely related to these \textit{forward iterations},
and in fact a key tool to the analysis of the ergodic property, is the
sequence of \textit{backward iterations}
\begin{equation}\label{3.6}
 \tilde{M}_{n} :=  F_{1:n}(M_{0}),\qquad  n\ge 0.
\end{equation}
The connection is established by the identity
\begin{equation}\label{3.7}
 \pi(y) {\mathbf P}(M_{n}\in\cdot|Y_0=y) = \pi(z) \tilde{\mathbf P}(\tilde{M}_n \in \cdot|\tilde{Y}_0=z)
\end{equation}
for all $n\ge 0$. Put also $M_{n}^{u} := F_{n:1}(u)$ and
$\tilde{M}_n^{u}:= F_{1:n}(u)$ for $u \in {\mathbf M}$ and note that
\begin{eqnarray}\label{3.8}
&&  \int_{z \in {\mathcal  Y}} \int_{y \in {\mathcal  Y}}   {\mathbf P}
\bigl((M_{n}^{u},\tilde{M}_n^u)_{n\ge 0}\in\cdot|Y_0=y,
\tilde{Y}_0=z\bigr) \pi(dy) \pi(dz) \nonumber
\\[-8pt]
\\[-8pt]
\nonumber &&\qquad = \int_{z \in {\mathcal  Y}} \int_{y \in {\mathcal
Y}} {\mathbf P}\bigl((M_{n},\tilde{M}_n)_{n\ge 0}
 \in\cdot|Y_0=y, \tilde{Y}_0=z\bigr) \pi(dy) \pi(dz).
\end{eqnarray}
Note that in (\ref{3.8}), the probability ${\mathbf P}$ denotes a joint
probability.

$\{Y_n,n \geq 0\}$ is called \textit{Harris recurrent} if there exist
a set $A \in {\mathcal  A}$, a probability measure $\Gamma$
concentrated
on~$A$ and an $\varepsilon$ with $0 < \varepsilon < 1$ such that
$P_y(Y_n \in A  \mbox{ i.o.})=1$ for all $y \in {\mathcal  Y}$ and,
furthermore, there exists $n$ such that $P^n(y,A') \geq  \varepsilon
\Gamma(A')$ for all $y \in A$ and all $A' \in {\mathcal  A}$.

A central question for an MIRFS $(M_n)_{n \ge 0}$ is under which
conditions it stabilizes, that is, converges to a stationary
distribution $\Pi$. The next theorem summarizes the results regarding
this question.

\begin{theorem}\label{T1}
Let $\{Y_n, n \geq 0\}$ be an aperiodic, irreducible and Harris
recurrent Markov chain, and let $(M_n)_{n \ge 0}$ be an MIRFS of
Lipschitz functions. Suppose the initial distribution of $Y_0$ is
$\pi$, and
\begin{equation}\label{3.10}
 {\mathbf E} \log l(F_{1}) < 0\quad\mbox{and}\quad {\mathbf E}
 \log^{+}d\bigl(F_{1}(u_{0}),u_{0}\bigr)<\infty
\end{equation}
for some $u_{0} \in {\mathbf M}$. Then the following assertions hold:
\begin{longlist}[(iii)]
\item[(i)] $\tilde{M}_{n}$ converges a.s. to a random element
$\tilde{M}_{\infty}$
 which does not depend on the initial distribution.

\item[(ii)] $M_n$ converges in distribution to $\tilde{M}_{\infty}$
under ${\mathbf P}$.

\item[(iii)] Define $\Pi$ as the stationary distribution of
$(\tilde{Y}_\infty, \tilde{M}_{\infty})$. Then $\Pi$ is the unique
stationary probability of the Markov chain $\{(Y_n, M_n), n \geq 0\}$.

\item[(iv)] $(M_n)_{n \ge 0}$ is ergodic under ${\mathbf P}_{\Pi}$,
that is, for any $u \in {\mathbf M}$,
\begin{equation}\label{3.11}
\frac{1}{n} \sum_{k=1}^n g(M_k) \longrightarrow {\mathbf E}_{\Pi}
(g(\tilde{M}_{\infty})),\qquad  {\mathbf P}_{\Pi}\mbox{-a.s.}
\end{equation}
for all bounded continuous real-valued functions $g$ on ${\mathbf M}$.
\end{longlist}
\end{theorem}

We remark that Elton \cite{E90} showed in the situation of a stationary
sequence $(F_{n})_{n\ge 1}$ that Theorem~\ref{T1} holds whenever
${\mathbf E} \log^{+}l(F_{1})$ and ${\mathbf
E}\log^{+}d(F_{1}(u_{0}),u_{0})$ are both finite for some (and then
all) $u_{0} \in {\mathbf M}$ and the Lyapunov exponent $ \gamma  :=
\lim_{n\to\infty} n^{-1}\log l(F_{n:1})$, which exists by Kingman's
subadditive ergodic theorem, is a.s. negative. Since the initial
distribution of $Y_0$ is the stationary distribution $\pi$, the Markov
chain $Y_n$ is a stationary sequence, and hence, $M_n$ is a sequence of
iterated random functions generated by stationary sequences. Here, we
impose the Harris recurrent condition so that the invariant measure
$\pi$ exists, and we are able to characterize $\tilde{M}_{\infty}$ in a
Markovian setting. Since the proof is similar to that in~\cite{BEH89}, it is
omitted.

\section{Central limit theorem and Edgeworth expansion for distributions
of a Markovian iterated random functions system}\label{s4}
Consider the Markovian
iterated random functions system $\{(Y_n,M_n),n\geq 0\}$ defined in
(\ref{3.1}). Abuse the notation a little bit and let $g$ be an ${\mathbf
R}^p$-valued function on ${\mathbf M}$. In this section we study
the central limit theorem and Edgeworth expansion of the sum
$S_n=\sum_{k=1}^n g(M_k)$ and ${\mathbf g}(n^{-1}S_n)$ for a smooth
function ${\mathbf g}\dvtx {\mathbf R}^p \to {\mathbf R}^q$. Let
$w\dvtx  {\mathcal Y} \rightarrow [1,\infty)$ be a measurable function,
and let ${\mathbf B}$ be the Banach space of measurable functions
$h\dvtx {\mathcal Y} \rightarrow C $ ($:=$ the set of complex numbers)
with $\|h\|_w := \sup_y|h(y)|/w(y) <\infty$. Assume further that
$\{Y_n,n\geq 0\}$ has a stationary distribution $\pi$ with $\int
w(y)\pi(dy) <\infty,$ and
\begin{eqnarray}
\hspace*{10mm} &\displaystyle{\lim_{n \rightarrow \infty} \sup_y \biggl\{\bigg|E[h(Y_n)|Y_0=y] - \int
  h(z)\pi(dz)\bigg|\Big/w(y)\dvtx y \in {\mathcal  Y}, |h| \leq w\biggr\} = 0,}& \label{4.1}
\\
&\displaystyle{\sup_y\{E[w(Y_p)|Y_0=y]/w(y)\} < \infty,}& \label{4.2}
\end{eqnarray}
for some $p \geq 1$. Condition (\ref{4.1}) says that the chain is
$w$-uniformly ergodic, which implies that there exist $\gamma > 0 $ and
$ 0<\rho<1 $ such that, for all $h\in {\mathbf B} $  and $ n \geq 1$,
\begin{equation}\label{4.3}
\sup_y \bigg|E[h(Y_n)|Y_0=y] - \int h(z)\pi(dz)\bigg|\Big/w(y) \leq \gamma
\rho^n\|h\|_w,
\end{equation}
(cf. pages 382--383 and Theorem~16.0.1 of~\cite{MT93}). We remark that,
for $w=1,$ condition (\ref{4.1}) is the classical uniform ergodicity
condition for $\{Y_n, n \geq 0 \}$.

The following assumption will be assumed throughout this section.

\renewcommand{\theass}{K}
\begin{ass}\label{assK}

K1.\label{assK1} Let $\{Y_n,n \geq 0\}$ be an aperiodic, irreducible
Markov chain satisfying conditions (\ref{4.1})--(\ref{4.2}).
Furthermore, we assume the initial distribution of $Y_0$ is $\pi$.

K2.\label{assK2} The MIRFS $(M_n)_{n \ge 0}$ has the weighted mean
contraction property, that is, there exists a $p \geq 1$ such that
\[
 \sup_y  \biggl\{{\mathbf E} \biggl( \log \frac{L_{p}w(Y_p)}{w(y)}\Big| Y_0=y \biggr)  \biggr\} < 0.
\]

K3.\label{assK3} There exists $u_0 \in {\mathbf M}$ for
which
\[
{\mathbf E} d^2\bigl(F_{1}(u_{0}),u_{0}\bigr)<
\infty\quad\mbox{and}\quad \sup_{y} \biggl\{{\mathbf E} \biggl( \frac{
L_1 w(Y_1)}{w(y)}\Big|Y_0=y \biggr) \biggr\} < \infty.
\]
\end{ass}

\begin{rem}\label{r1}
(a) Assumption~\hyperref[assK1]{K1} is a condition for the underlying
Markov chain $\{Y_n,n \geq 0\}$ which is general enough to include
several practical used models studied in Section~\textup{\ref{s6}}.
Assumption~\hyperref[assK2]{K2} is a weighted mean contraction
condition which is different from the standard mean contraction
condition \mbox{${\mathbf E} \log L_1 < 0$} used in Theorem~\ref{T1}.
Assumption \hyperref[assK3]{K3} is a weighted moment condition. Note that under
Assumptions \hyperref[assK1]{K1}--\hyperref[assK3]{K3}, and the extra assumption
that $\{(Y_n,M_n), n \geq 0\}$ is an irreducible, aperiodic and Harris
recurrent Markov chain, Theorems 13.0.1~and~17.0.1(i) of~\cite{MT93}
imply that Theorem~\ref{T1} still holds. Furthermore, we will prove the
central limit theorem and Edgeworth expansion for the distributions of
a Markovian iterated random functions system in Theorem~2.

(b) To have better understanding of Assumption~\ref{assK}, we consider
a simple state space model. Given $p \geq 1$ as in
Assumption \hyperref[assK2]{K2}, and $|\alpha| < 1$, let $Y_n = \alpha Y_{n-1} +
\varepsilon_n, \xi_n=\beta_{Y_n} \xi_{n-1} + \eta_n$, where
$\varepsilon_n$ are i.i.d. random variables with $E|\varepsilon_1|=c <
\infty$, and $\eta_n$ are i.i.d. random variables with $E|\eta_1|<
\infty$. Further, we assume both $\varepsilon_1$ and $\eta_1$ have
positive probability density function with respect to Lebesgue
measure, and that they are mutually independent. Denote
$b=(1-|\alpha|^p)/(1-|\alpha|)$ and $a=1/(bc + 1) < 1$, and assume
$|\beta_y| < a^{1/p}<1$ for all $y \in {\mathcal  Y}$. It is known that
$w(y) = |y| +1$ (cf. pages 380~and~383 of~\cite{MT93}). Let
$d(u,v)=|u-v|.$ It is easy to see that Assumption \hyperref[assK1]{K1}
and the first part of Assumption \hyperref[assK3]{K3} hold. To check
Assumption~\hyperref[assK2]{K2}, we have
\begin{eqnarray}\label{4.3a}
\nonumber && \sup_y  \biggl\{{\mathbf E} \biggl(\log \frac{L_p
w(Y_p)}{w(y)}\Big|Y_0=y \biggr) \biggr\}
\\
&&\qquad = \sup_y  \biggl\{{\mathbf E} \biggl(\log
\frac{|\beta_{Y_p}\cdots \beta_{Y_1}|(|\alpha^p y + \sum_{k=0}^{p-1}
\alpha^k
\varepsilon_{p-k}|+1)}{|y|+1}\Big|Y_0=y \biggr) \biggr\}
\nonumber
\\[-8pt]
\\[-8pt]
\nonumber &&\qquad < \log   \sup_{y}  \biggl\{  \frac{a ( |\alpha^p y|
+ E|\sum_{k=0}^{p-1} \alpha^k \varepsilon_{p-k}|+1)}{|y|+1}  \biggr\}
\\
\nonumber &&\qquad = \log \sup_{y} \biggl\{ \frac{a(|\alpha^p y| +  bc
+ 1)}{|y|+1} \biggr\} = 0.
\end{eqnarray}
By using the same argument, we have the second part of
Assumption \hyperref[assK3]{K3}. When $\varepsilon_n$ are i.i.d. $N(0,1)$,
$\eta_n$ are i.i.d. $N(0,1)$, and they are mutually independent. Then
$a=\sqrt{2\pi}/(2b + \sqrt{2\pi}\,)<1$.

Recall that $\Pi$ is defined in Theorem~\ref{T1}(iii) and denote
$Q(B):=\Pi({\mathcal  Y} \times B)$ for all $B \in {\mathcal
B}({\mathbf M})$. Let $g \in {\mathcal  L}_{0}^{2}(Q)$ be a square
integrable function taking values in ${\mathbf R}^p$ with mean $\bolds{0}$,
that is, $g=(g_1,\ldots,g_p)$ with each $g_k$ a real-valued function on
${\mathbf M}$, and
\begin{equation}\label{4.4}
 \int_{\mathbf M} g_k(u) Q(du) = 0,\qquad  \|g_k\|_{2}^{2} = \int_{\mathbf M} g_k^{2}(u) Q(du) < \infty,
\end{equation}
for $k=1,\ldots,p$. Consider the sequence
\begin{equation}\label{4.5}
S_n= S_{n}(g) = g(M_{1})+\cdots + g(M_{n}),\qquad  n\geq 1,
\end{equation}
which may be viewed as a Markov random walk on the Markov chain
$\{(Y_n,M_n), n \geq 0\}$.

Note that there are two special properties of the Markov chain induced
by the Markovian iterated random functions system
(\ref{2.4})--(\ref{2.7}). First, the hypothesis that the transition
probability possesses a density leads to a classical situation in the
context of the so-called ``Doeblin condition'' for Markov chains.
Second, a positivity hypothesis on ${\mathbf M}$ defined in (\ref{2.4})
in the support of the Markov chain leads to contraction properties, on
which basis we will develop the spectral theory. The reader is referred
to \cite{HH01} for a general account of the
perturbation theory of Markovian operators. We need the following
notation first.
\end{rem}

\begin{definition}\label{D2}
Let $w\dvtx {\mathcal  Y} \to [1,\infty)$ be a weight function. For
any measurable function $\varphi\dvtx  {\mathcal  Y} \times {\mathbf M}
\rightarrow [1,\infty)$, given $u_0 \in {\mathbf M}$, define
\[
 \|\varphi\|_w := \sup_{ y \in {\mathcal  Y}, u \in {\mathbf M}}
  \frac{|\varphi(y,u)|}{w(y)}
\]
and
\[
\|\varphi\|_h := \sup_{y \in {\mathcal  Y}, u,v: 0 < d(u,v) \leq 1}
\frac{|\varphi(y,u) - \varphi(y,v)|}{(w(y)\,d(u,v))^{\delta}},
\]
for $0 < \delta <1$. We define ${\mathcal  H}$ as the set of $\varphi$
on ${\mathcal  Y} \times {\mathbf M}$ for which $\|\varphi\|_{wh} :=
\|\varphi\|_w + \|\varphi\|_h$ is finite, where $wh$ represents a
combination of the weighted variation norm and the bounded weighted
H\"{o}lder norm.
\end{definition}

Let $\nu$ be an initial distribution of $(Y_0,M_0)$ and let ${\mathbf
E}_\nu$ denote expectation under the initial distribution $\nu$ on
$(Y_0,M_0)$. For $\varphi \in {\mathcal  H}$, $g \in {\mathcal
L}^{2}(Q)$, $y \in {\mathcal  Y}$, $u \in {\mathbf M}$ and $p \times
1$ vectors $\alpha=(\alpha_1,\ldots,\alpha_p)' \in {\mathbf R}^p$,
define linear operators ${\mathbf T}_{\alpha}$, ${\mathbf T}$,
$\nu_{\alpha}$ and ${\mathbf Q}$ on the space ${\mathcal  H}$ as
\begin{eqnarray}
({\mathbf T}_{\alpha} \varphi)(y,u) &=& {\mathbf E}\bigl\{e^{i \alpha'
g(M_1)}     \varphi(Y_1, M_1)|Y_0=y,M_0= u \bigr\}, \label{4.11}
\\
({\mathbf T} \varphi)(y,u) &=& {\mathbf E} \{\varphi(Y_1, M_1)|Y_0=y,M_0= u \}, \label{4.12}
\\
\nu_{\alpha}\varphi &=& {\mathbf E}_\nu \bigl\{e^{i\alpha' \varphi(u)}
\varphi(Y_0,u) \bigr\},\qquad   {\mathbf Q}\varphi = {\mathbf E}_\Pi
\{\varphi(Y_0,u)\}. \label{4.13}
\end{eqnarray}
In the case of a $w$-uniformly ergodic Markov chain, Fuh and Lai
\cite{FL01} have shown that there exists a sufficiently small $\delta
> 0$ such that, for $|\alpha| \leq \delta$, ${\mathcal  H}={\mathcal  H}_1(\alpha)
\oplus {\mathcal  H}_2(\alpha)$ and
\begin{equation}
{\mathbf T}_\alpha {\mathbf Q}_\alpha \varphi = \lambda(\alpha)
{\mathbf Q}_\alpha \varphi\qquad\mbox{for all }  \varphi \in {\mathcal
H}, \label{4.14}
\end{equation}
where ${\mathcal  H}_1(\alpha)$ is a one-dimensional subspace of
${\mathcal  H}$, $\lambda(\alpha)$ is the eigenvalue of ${\mathbf
T}_\alpha$ with corresponding eigenspace ${\mathcal  H}_1(\alpha)$
and~${\mathbf Q}_{\alpha}$ is the parallel projection of ${\mathcal H}$
onto the subspace ${\mathcal  H}_1(\alpha)$ in the direction of
${\mathcal H}_2(\alpha).$ Extension of their argument to the weight
functions $w$ and $l$ defined in Definition \ref{D2} is given in the
\hyperref[app]{Appendix}, which also proves the following lemmas.

\begin{lemma}\label{L1}
Let $\{(Y_n,M_n),n\geq 0\}$ be the MIRFS of Lipschitz functions defined
in \textup{(\ref{2.1})} and satisfying Assumption~\textup{\ref{assK}}. Assume $g
\in {\mathcal L}^{r}(Q)$ for some $r > 2$. Then ${\mathbf T}$~and~${\mathbf Q}$
are bounded linear operators on the Banach space
${\mathcal  H}$ with norm $\| \cdot \|_{wh}$, and satisfy
\begin{equation}\label{4.15}
\| {\mathbf T}^n- {\mathbf Q} \|_{wh} =\sup_{\varphi \in {\mathcal H},
\|\varphi\|_{wh} \leq 1} \| {\mathbf T}^n \varphi - {\mathbf Q} \varphi
\|_{wh} < \gamma_* \rho_*^n,
\end{equation}
for some $\gamma_* > 0$ and $0 < \rho_* <1$.
\end{lemma}

By using an argument similar to Proposition 1 of~\cite{F04b}, we have
the following:

\begin{lemma}\label{L2}
Let $\{(Y_n,M_n),n\geq 0\}$ be the MIRFS defined in \textup{(\ref{2.1})}
satisfying Assumption \textup{\ref{assK}}, such that the induced Markov
chain $\{(Y_n,M_n), n\geq 0\}$ with transition probability kernel
\textup{(\ref{3.2})} is irreducible, aperiodic and Harris recurrent. Assume
$g \in {\mathcal L}^{r}(Q)$ for some $r > 2$. Then there exists
$\delta > 0 $ such that, for $\alpha \in {\mathbf R}^p$ with $|\alpha|
< \delta $, and for $\varphi \in {\mathcal  H},$
\begin{eqnarray}\label{4.16}
 {\mathbf E}_\nu \bigl\{e^{i \alpha' g(M_n)} \varphi(Y_n,M_n)\bigr\}
&=& \nu_\alpha {\mathbf T}_\alpha^n \varphi = \nu_\alpha {\mathbf
T}_\alpha^n\{{\mathbf Q}_\alpha +(I-{\mathbf Q}_\alpha)\}\varphi
\nonumber
\\[-8pt]
\\[-8pt]
\nonumber  &=& \lambda^n(\alpha)\nu_\alpha {\mathbf Q}_\alpha \varphi +
\nu_\alpha {\mathbf Q}_\alpha^n(I- {\mathbf Q}_\alpha)\varphi,
\end{eqnarray}
and:
\begin{longlist}[(iii)]
\item[(i)] $\lambda(\alpha)$ is the unique eigenvalue of the maximal
modulus of ${\mathbf T}_\alpha$;

\item[(ii)] ${\mathbf Q}_\alpha$ is a rank-one projection;

\item[(iii)] the mappings $\lambda(\alpha), {\mathbf Q}_\alpha$ and
$I-{\mathbf Q}_\alpha$ are analytic;

\item[(iv)] $|\lambda(\alpha)| > \frac{2+\rho_*}{3}$ and for each $k
\in N,$ the set of positive integers, there exists $c > 0$ such that,
for each $n \in N$ and $j_1,\ldots,j_p$ with $j_1+\cdots+j_p=k$,
\[
\bigg\|\frac{\partial^k}{\partial \alpha_1^{j_1} \cdots \partial
\alpha_p^{j_p}} (I-{\mathbf Q}_\alpha)^n\bigg\|_{wh} \leq
c\biggl(\frac{1+2\rho_*}{3}\biggr)^n;
\]

\item[(v)] denote $g=(g_1,\ldots,g_p)$, and let $\gamma_j :=
\lim_{n\rightarrow \infty} (1/n){\mathbf E}_{yu} \log\|g_j(M_n) \|$,
the upper Lyapunov exponent; it follows that
\begin{equation}\label{4.17}
 \gamma_j = \frac{\partial \lambda(\alpha)}{\partial \alpha_j}\bigg|_{\alpha =0}
= \int {\mathbf E}_{yu}g_j(M_1) \Pi(dy\times du).
\end{equation}
\end{longlist}
\end{lemma}

Note that in Lemma \ref{L2} we need the extra assumption that the
induced Markov chain $\{(Y_n,M_n),n\geq 0\}$ with transition
probability kernel ({\ref{3.2}) is irreducible, aperiodic and Harris
recurrent. In Section~\ref{s5} we will show that this condition is
satisfied for the Markov chain induced by the Markovian iterated random
functions system (\ref{2.4})--(\ref{2.7}).

For given $S_n=\sum_{k=1}^n g(M_k)$ of the MIRFS $\{(Y_n,M_n), n\geq
0\}$, in this section we will obtain Edgeworth expansions for the
standardized distribution of $S_n$ via the representation (\ref{4.16})
of the characteristic function ${\mathbf E}(e^{i\alpha'
g(M_n)}|Y_0=y,M_0=0)$. Note that Lemma \ref{L1} implies that
$\{(Y_n,M_n), n\geq 0 \}$ is geometrically mixing in the sense that
there exist $r_1 > 0$ and $0 <\gamma_1<1$ such that, for all $y \in
{\mathcal  Y},u \in {\mathbf M}, k \geq 0 $ and $n \geq 1 $ and for all
real-valued measurable functions $\varphi_1, \varphi_2$ with
$\|\varphi_1^2\|_{wh} < \infty $ and $\|\varphi_2^2\|_{wh} < \infty$,
\begin{eqnarray}\label{4.18}
\nonumber && \|{\mathbf E}\{\varphi_1(Y_k,M_k)\varphi_2(Y_{k+n}, M_{k+n})|Y_0=y,M_0=u\}
\\
&&\phantom{\|}{}   - \{{\mathbf E}
\varphi_1(Y_k,M_k)|Y_0=y,M_0=u\}
\\
\nonumber &&\hspace*{13mm} \phantom{\|}{} {}\times \{{\mathbf E}
\varphi_2(Y_{k+n},M_{k+n}|Y_0=y,M_0=u)\}\|_{wh} \leq r_1 \gamma_1^n.
\end{eqnarray}
Let $\tilde{\varphi}_1, \tilde{\varphi}_2$ be real-valued measurable
functions on $({\mathcal  Y}\times {\mathbf M}) \times ({\mathcal  Y}
\times {\mathbf M})$. Denote $\varphi_1(z,v) ={\mathbf
E}\{\tilde{\varphi}_1((z,v),(Y_1,M_1))|Y_0=z,M_0=v)\}$, and note that
\begin{eqnarray*}
&& {\mathbf
E}\bigl\{\tilde{\varphi}_1\bigl((Y_k,M_k),(Y_{k+1},M_{k+1})\bigr)|Y_0=y,M_0=u\bigr\}
\\
&&\qquad ={\mathbf E}\{\varphi_1(Y_k,M_k)|Y_0=y,M_0=u\}.
\end{eqnarray*}
The same proof as that of Theorem~16.1.5 of
\cite{MT93} can be used to show that there exist $r_1>0$ and
$0<\gamma_1 <1$ such that, for all $y \in {\mathcal  Y}, u \in {\mathbf
M}, k\geq 0$ and $ n\geq 1$ and for all measurable $\tilde{\varphi}_1,
\tilde{\varphi}_2$ with $\|\sup_{z,v}\tilde{\varphi}_1^2((y,u),
(z,v))\|_{wh} < \infty $ and $\|\sup_{z,v}
\tilde{\varphi}_2^2((y,u),(z,v))\|_{wh} < \infty,$
\begin{eqnarray}\label{4.19}
\nonumber \hspace*{9mm} && \big\|{\mathbf
E}\bigl\{\tilde{\varphi}_1\bigl((Y_k,M_k),(Y_{k+1},M_{k+1})\bigr)
\\
\nonumber &&\phantom{\big\|{\mathbf E}\bigl\{}
{}\times \tilde{\varphi}_2
\bigl((Y_{k+n},M_{k+n}),(Y_{k+n+1},M_{k+n+1})\bigr)|Y_0=y,M_0=u\bigr\}
\\
&&\phantom{\big\|}{}   - {\mathbf
E}\bigl\{\varphi(Y_k,M_k)|Y_0=y,M_0=u\bigr\}{\mathbf E}\{
\varphi_2(Y_{k+n},M_{k+n})|Y_0=y,M_0=u\}\big\|_{wh}
\hspace*{-12pt}\\
\nonumber &&\qquad \leq r_1 \gamma_1^{n-1}.
\end{eqnarray}

To establish Edgeworth expansion for a Markovian iterated random
functions system, we shall make use of (\ref{4.19}) in conjunction with
the following extension of Cram\'{e}r (strongly nonlattice)
condition:
\begin{equation}\label{4.20}
\inf_{|v| > \alpha} |1- E_{\pi} \{\exp(i v' S_1(g) \}| > 0\qquad
\mbox{for all } \alpha > 0.
\end{equation}
In addition, we also assume  the conditional Cram\'{e}r (strongly
nonlattice) condition ((2.5) on page 216 in \cite{GH83}):
There exists $\delta > 0$ such that, for all
$m,n=1,2,\ldots,$ $\delta^{-1} <m< n$, and all $\alpha \in {\mathbf
R}^p$ with $|\alpha| \geq \delta$,
\begin{eqnarray}\label{4.21}
\nonumber \hspace*{8mm} && E_{\pi} \big| E\bigl\{\exp\bigl(i\alpha' \bigl(g(M_{n-m}) +\cdots+ g(M_{n+m})\bigr)\bigr)
\\
&&\hspace*{10.5mm}
  \big|(Y_{n-m},M_{n-m}),\ldots,(Y_{n-1},M_{n-1}),
\\
\nonumber &&\hspace*{11.5mm}  (Y_{n+1},M_{n+1}),\ldots,(Y_{n+m},M_{n+m}),(Y_{n+m+1},M_{n+m+1})\bigr\}\big|
 \leq e^{-\delta}.
\end{eqnarray}

Let
\begin{equation}\label{4.22}
\gamma = \int {\mathbf E}_{yu} g(M_1) \Pi(dy \times du)   \bigl( =
\lambda'(0)\bigr),
\end{equation}
and denote by $V = (\partial^2\lambda(\alpha)/\partial\alpha_i\,\partial
\alpha_j|_{\alpha=0})_{1\leq i,j\leq p}$  the Hessian matrix of
$\lambda$ at 0. By Lemma \ref{L2},
\begin{equation}\label{4.23}
 \lim_{n\rightarrow\infty}n^{-1}{\mathbf E}_\nu\bigl\{\bigl(g(M_n)-n\gamma\bigr)
 \bigl(g(M_n)-n\gamma\bigr)'\bigr\}= V.
\end{equation}

Let $\psi_n(\alpha)={\mathbf E}_\nu(e^{i\alpha' g(M_n)})$. Then by
Lemma \ref{L2} and the fact that $\nu_\alpha {\mathbf Q}_\alpha h_1$
has continuous partial derivatives of order $r-2$ in some neighborhood
of $\alpha = 0 $, we have the Taylor series expansion of
$\psi_n(\alpha/\sqrt{n}\,)$ for $|\alpha/\sqrt{n}| \leq \varepsilon $
(some sufficiently small positive number):
\begin{equation}\label{4.24}
\psi_n\bigl(\alpha/\sqrt{n}\,\bigr)
\Biggl\{1+\sum_{j=1}^{r-2}n^{-j/2}\tilde{\pi}_j(i\alpha)\Biggr\}
e^{-\alpha'V\alpha/2}+ o\bigl(n^{-(r-2)/2}\bigr),
\end{equation}
where $\tilde{\pi}_j(i\alpha)$ is a polynomial in $i\alpha$ of degree
$3j$ whose coefficients are smooth functions of the partial
derivatives of $\lambda(\alpha) $ at $\alpha = 0 $ up to the order
$j+2$ and those of $\nu_\alpha {\mathbf Q}_\alpha h_1$ at $\alpha = 0$
up to the order $j$. Letting $D$ denote the $p \times 1$ vector whose
$j$th component is the partial differentiation operator $D_j$ with
respect to the $j$th coordinate, define the differential operator
$\tilde{\pi}_j(-D)$. As in the case of sums of i.i.d. zero-mean random
vectors (cf. \cite{BR76}), we obtain an Edgeworth expansion for the
``formal density'' of the distribution of $g(M_n)$ by replacing the
$\tilde{\pi}_j(i\alpha)$ and $e^{-\alpha'V\alpha/2}$ in (\ref{4.24}) by
$\tilde{\pi}_j(-D)$ and $\phi_V(y)$, respectively, where $\phi_V$ is
the density function of the $q$-variate normal distribution with mean 0
and covariance matrix~$V$. Throughout the sequel we let ${\mathbf
P}_{\nu}$ denote the probability measure under which $(Y_0,M_0)$ has
initial distribution $\nu$.

\begin{theorem}\label{T2}
Let $\{(Y_n,M_n),n\geq 0\}$ be the MIRFS defined in
\textup{(\ref{2.1})} satisfying Assumption~\textup{\ref{assK}}, such
that the induced Markov chain $\{(Y_n,M_n),n\geq 0\}$, with transition
probability kernel~\textup{(\ref{3.2})}, is irreducible, aperiodic and
Harris recurrent. Assuming $g \in {\mathcal  L}^{r}(Q)$ for some $r >
2$, \textup{(\ref{4.20})} and~\textup{(\ref{4.21})} hold. Let
$\phi_{j,V}=\tilde{\pi}_j(-D)\phi_V $  for $j = 1,\ldots,r-2$. For $0<a
\leq 1$ and $c > 0 $, let ${\mathcal  B}_{a,c}$ be the class of all
Borel subsets $B$ of $\mathbf{ R}^p$ such that $\int_{(\partial
B)^{\varepsilon}}\phi_V(y)\,dy \leq c\varepsilon^{a}$ for every
$\varepsilon > 0$, where $\partial B$ denotes the boundary of $B$ and
$(\partial B)^\varepsilon$ denotes its $\varepsilon$-neighborhood. Then
\begin{eqnarray}\label{4.25}
\hspace*{6mm} &&  \sup_{B \in{\mathcal  B}_{a,c}}\Bigg|{\mathbf P}_\nu\bigl\{(S_n-n\gamma)/\sqrt{n} \in B\bigr\}-
\int_B \Biggl\{\phi_V(y)+\sum_{j=1}^{r-2}n^{-j/2}\phi_{j,V}(y)\Biggr\}\,dy\Bigg|
\nonumber
\\[-2pt]
\\[-14pt]
\nonumber
&&\qquad =o\bigl(n^{-(r-2)/2}\bigr).
\end{eqnarray}
\end{theorem}

A proof of Theorem~2 is given in the \hyperref[app]{Appendix}.

Note that under weaker moment conditions, and an alternative condition
of (\ref{4.20}) and (\ref{4.21}) (see Condition 1 of~\cite{L93}) Lahiri
\cite{L93} proved the asymptotic expansions for sums of weakly
dependent random vectors.

Letting  $r=2$ in Theorem~\ref{T2}, we have the following:

\begin{corollary}\label{C1}
With the same notation and assumptions as in~Theorem~\textup{\ref{T2}},
then
\[
 \frac{1}{\sqrt{n}}  (S_n  - n \gamma  )  \longrightarrow
N(0,\Sigma)\qquad\mbox{in distribution},
\]
where the variance--covariance matrix
\begin{equation}\label{4.26}
\Sigma =  \biggl( \frac{ \partial^2 \lambda(\alpha)}{ \partial
\alpha_i\,
\partial \alpha_j}\bigg|_{\alpha=0 }
 \biggr)_{i,j=1,\ldots,p}.
\end{equation}
\end{corollary}

In statistical applications one often works with ${\mathbf
g}(n^{-1}S_n)$ instead of $S_n=\sum_{k=1}^n g(M_k)$, where ${\mathbf
g}\dvtx  {\mathbf R}^p \rightarrow {\mathbf R}^q$ is sufficiently
smooth in some neighborhood of the mean $\gamma:=(
\gamma_1,\ldots,\gamma_p)$. Denote ${\mathbf g} = ({\mathbf
g}_1,\ldots,{\mathbf g}_q)$ with each ${\mathbf g}_i$, $1 \leq i \leq
q$, a real-valued function on ${\mathbf R}^p$. For the case of a sum of
i.i.d. random variables, Bhattacharya and Ghosh \cite{BG78} made use of
the Edgeworth expansion of the distribution of $(S_n-n\gamma)/\sqrt{n}$
to derive an Edgeworth expansion of the distribution of
$\sqrt{n}\{{\mathbf g}(n^{-1}S_n)- {\mathbf g}(\gamma)\}$. Making use
of Theorem~\ref{T2} and a straightforward extension of their argument,
we can generalize their result to the case where $S_n$ is the partial
sum of a Markovian iterated random functions system.

\begin{theorem}\label{T3}
Under the same assumptions as in~Theorem~\textup{\ref{T2}}, suppose
that ${\mathbf g} \dvtx {\mathbf R}^p \rightarrow {\mathbf R}^q$ has
continuous partial derivatives of order $r$ in some neighborhood of
$\gamma$. Let $J_{\mathbf g}=(D_j {\mathbf g}_i(\gamma))_{1\leq i \leq
q,1\leq j \leq p}$ be the $q\times p$ Jacobian matrix and let
$V({\mathbf g}) = J_{\mathbf g} V J_{\mathbf g}'$. Then
\begin{eqnarray}\label{4.27}
\nonumber && \sup_{B\in {\mathcal  B}_{a,c}}\Bigg|{\mathbf
P}_\nu\bigl\{\sqrt{n}\bigl({\mathbf g}(n^{-1}S_n)- {\mathbf
g}(\gamma)\bigr) \in B \bigr\}
\\
&& \hspace*{10mm}
{}- \int_B \Biggl\{\phi_{V({\mathbf
g})}(y)+ \sum_{j=1}^{r-2}n^{-j/2}\phi_{j,V,{\mathbf
g}}(y)\Biggr\}\,dy\Bigg|
\\
\nonumber &&\qquad = o\bigl(n^{-(r-2)/2}\bigr),
\end{eqnarray}
where $\phi_{j,V,{\mathbf g}} = \tilde{\pi}_{j,{\mathbf g}}(-D)\phi_V $
and $\tilde{\pi}_{j,{\mathbf g}}(y)$ is a polynomial in $y (\in
{\mathbf R}^p)$ whose coefficients are smooth functions of the partial
derivatives of $\lambda(\alpha)$ at $\alpha =0$  up to order $j+2$ and
those of $\nu_\alpha {\mathbf Q}_\alpha h_1$ at $\alpha = 0$ up to
order $j$ together with those of ${\mathbf g}$ at $\mu$ up to order $j
+1$.
\end{theorem}

In the next theorem we consider $p=1$.

\begin{theorem}\label{T4}
Under the same assumptions as in~Theorem~\textup{\ref{T2}}, assume $g
\in {\mathcal  L}^{r}(Q)$ for some $r > 2$. Then
\begin{equation}\label{4.28}
\hspace*{10mm} \frac{1- {\mathbf P}_\nu\{(S_n-n\gamma)/\sqrt{n} \leq t\}}{1 - \Phi(t)}
= \exp\bigl(t^3/\sqrt{n}\,\bigr) \varphi\bigl(t/\sqrt{n}\,\bigr)
\biggl(1 + O\biggl(\frac{t}{\sqrt{n}}\biggr) \biggr)
\end{equation}
and
\begin{equation}\label{4.29}
\hspace*{11mm}  \frac{{\mathbf P}_\nu\{(S_n-n\gamma)/\sqrt{n} \leq -t\}}{\Phi(-t)} =
\exp\bigl(-t^3/\sqrt{n}\,\bigr) \varphi\bigl(-t/\sqrt{n}\,\bigr)
\biggl(1 + O\biggl(\frac{t}{\sqrt{n}}\biggr) \biggr),
\end{equation}
where $\Phi(t)$ is the standard normal distribution, and $\varphi(t)$
is a power series which converges for $t$ sufficiently small in absolute
value.
\end{theorem}

Theorem~\ref{T4} states the moderate deviations results for
the distribution of an MIRFS, which will be used to prove Edgeworth
expansion for the MLE in Section~\ref{s5}. Since the proof is a
straightforward generalization of Theorem~6 in~\cite{N61}, it will not
be repeated here.

\section{Efficient likelihood estimation}\label{s5}
For a given state space model defined in~(\ref{2.1}) which involves
several parameters $\theta =( \theta_1,\ldots,\theta_q)$, the
estimation problem we consider in this section is the case of
estimating one of the parameters at a time; the other parameters play
the role of nuisance parameters. The true parameter is denoted by
$\theta_0$. Recall $p_n=p_n(\xi_0,\xi_1,\ldots,\xi_n;\theta)$ defined
as (\ref{2.3}). When $\partial \log p_n/\partial \theta$ exists, one
can seek solutions of the likelihood equations
\begin{equation}\label{5.1}
 \frac{\partial \log p_n}{\partial \theta} = 0.
\end{equation}

In the following, we denote $E^{\theta}_x$ as the expectation defined
under $P^{\theta}(\cdot,\cdot)$ in~(\ref{2.1}) with initial state
$X_0=x,$ and $E^\theta_{(x,s)}$ as the expectation defined under
$P^{\theta}(\cdot,\cdot)$ in~(\ref{2.1}) with initial state
$X_0=x,\xi_0=s$. The following conditions will be used throughout the
rest of this paper.

C1. For given $\theta \in \Theta$, the Markov chain $\{(X_n,\xi_n), n
\geq 0 \}$ defined in (\ref{2.1}) and~(\ref{2.2})
is aperiodic, irreducible, and satisfies (\ref{4.1})~and~(\ref{4.2}) with weight
function $w(\cdot)$. Assume $0< p_{\theta}(x,y)< \infty$ for all $x,y
\in {\mathcal  X}$, and $0< \sup_{x \in {\mathcal  X}}
f(s_1;\break \theta|x,s_0)< \infty,$ for all $s_0,s_1 \in {\mathbf R}^d$.
Denote $g_{\theta}(s_0,\xi_1)= \sup_{x_{0} \in {\mathcal  X}} \int
p_{\theta}(x_{0},x_1)\* f(\xi_1;\theta|x_1,s_{0})
  m(d x_1)$. Furthermore,
we assume that there exists $p \geq 1$ as in Assumption~\hyperref[assK2]{K2} such
that
\begin{eqnarray}
  \sup_{(x_0, s_0) \in {\mathcal  X} \times {\mathbf R}^d}
E^{\theta}_{(x_0,s_0)}  \biggl\{ \log  \biggl( g_{\theta}(s_0,\xi_1)^p
\frac{w(X_p,\xi_p)}{w(x_0,s_0)}   \biggr)  \biggr\}  &<& 0 \label{5.1a},
\\
\sup_{(x_0,s_0) \in {\mathcal  X} \times {\mathbf R}^d}
E^{\theta}_{(x_0,s_0)}  \biggl\{ g_{\theta}(s_0,\xi_1)
\frac{w(X_1,\xi_1)}{w(x_0,s_0)}  \biggr\}  &<& \infty. \label{5.1b}
\end{eqnarray}

C2. The true parameter $\theta_0$ is an interior point of $\Theta$. For
all $x\in {\mathcal  X}$, $s_0,s_1 \in {\mathbf R}^d$, $\theta \in
\Theta \subset {\mathbf R}^q$, and for $i,j,k =1,\ldots,q$, the partial
derivatives
\[
\frac{\partial {f(s_0;\theta|x)}}{\partial {\theta_i}},\qquad
\frac{\partial^2  {f(s_0;\theta|x)}} {\partial {\theta_i}\, \partial
{\theta_j}},\qquad  \frac{\partial^3 f(s_0;\theta|x)} {\partial
\theta_i\,\partial \theta_j\, \partial \theta_k}\qquad  \mbox{exist,}
\]
as well as the partial derivatives
\[
\frac{\partial {f(s_1;\theta|x,s_0)}}{\partial {\theta_i}},\qquad
\frac{\partial^2  {f(s_1;\theta|x,s_0)}}
  {\partial {\theta_i}\, \partial {\theta_j}},\qquad   \frac{\partial^3
f(s_1;\theta|x,s_0)} {\partial \theta_i\, \partial \theta_j\, \partial
\theta_k},
\]
and for all $x,y \in {\mathcal  X}$, $\theta \rightarrow p_\theta(x,y)$
and $\theta \rightarrow \pi_{\theta}(x)$ have twice continuous
derivatives in some neighborhood $N_\delta(\theta_0):=\{\theta\dvtx
|\theta-\theta_0| < \delta\}$ of $\theta_0$.

C3.
\[
\int_{\mathcal  X} \sup_{\theta \in N_\delta(\theta_0)} \bigg|
\frac{\partial \pi_\theta(x)}{\partial \theta_i}\bigg| m(dx) <
\infty,\qquad
 \int_{\mathcal  X} \sup_{\theta \in N_\delta(\theta_0)} \bigg|\frac{\partial^2  \pi_{\theta}(x)}
 {\partial \theta_i\, \partial \theta_j}\bigg| m(dx)<\infty,
\]
 and for all $x \in {\mathcal  X}$, $i,j =1,\ldots,q$,
\[
\int_{\mathcal  X} \sup_{\theta \in N_\delta(\theta_0)} \bigg|
\frac{\partial p_{\theta}(x,y)}
 {\partial \theta_i}\bigg| m(dy) < \infty,\qquad
 \int_{\mathcal  X} \sup_{\theta \in N_\delta(\theta_0)}
 \bigg|\frac{\partial^2  p_{\theta}(x,y)}{\partial \theta_i\, \partial \theta_j}\bigg| m(dy) < \infty.
\]

C4. For all $x \in {\mathcal  X}$, $s_0 \in {\mathbf R}^d$ and $\theta
\in \Theta$,
\begin{eqnarray*}
E_x^{\theta} \bigg| \frac{\partial f(\xi_0;\theta|x)} {\partial
\theta_i}\bigg| &<& \infty,\hspace*{18mm}      E_x^{\theta}
\bigg|\frac{\partial^2 f(\xi_0;\theta|x)} {\partial \theta_i\, \partial
\theta_j}\bigg| < \infty,
\\
E_{(x,s_0)}^{\theta} \bigg| \frac{\partial f(\xi_1;\theta|x,s_0)}
{\partial \theta_i}\bigg| &<& \infty,\qquad   E_{(x.s_0)}^{\theta}
\bigg|\frac{\partial^2 f(\xi_1;\theta|x,s_0)} {\partial \theta_i\,
\partial \theta_j}\bigg| < \infty.
\end{eqnarray*}
Furthermore, we assume that, for all $x \in {\mathcal  X}$, $s_0 \in
{\mathbf R}^d$ and uniformly for $\theta \in N_\delta(\theta_0)$,
\[
\bigg|\frac{\partial^3  \log f(\xi_0;\theta|x)} {\partial \theta_i\,
\partial \theta_j\, \partial \theta_k}\bigg| < H_{ijk}(x,\xi_0),\qquad
\bigg|\frac{\partial^3  \log f(\xi_1;\theta|x,s_0)} {\partial
\theta_i\, \partial \theta_j\, \partial \theta_k}\bigg| < G_{ijk}\bigl((x,s_0),\xi_1\bigr),
\]
where $H_{ijk}$ and $G_{ijk}$ are such that $E_x^{\theta_0}
H_{ijk}(x,\xi_0)< \infty$ and $E_{(x,s_0)}^{\theta_0}
G_{ijk}((x,s_0),\break \xi_1)< \infty$, for all $i,j,k =1,\ldots,q$ and for
all $x \in {\mathcal  X}, s_0 \in {\mathbf R}^d$.

C5.
\[
\sup_{x \in {\mathcal  X}} E_{x}^{\theta_0}  \biggl(\, \sup_{|\theta-
\theta_0| < \delta} \sup_{y,z \in {\mathcal  X}}
\frac{f(\xi_0;\theta|y)f(\xi_1;\theta|y,\xi_0)}
{f(\xi_0;\theta|z)f(\xi_1;\theta|z,\xi_0)}  \biggr)^2  < \infty.
\]

C6. The equality
\[
p_n(\xi_0,\xi_1,\ldots,\xi_n;\theta) = p_n(\xi_0,\xi_1,\ldots,\xi_n;\theta')
\]
holds $P$-almost surely, for all nonnegative $n$, if and only if
$\theta=\theta'$.

C7. For all $x,y \in {\mathcal  X}$, $\theta \rightarrow
p_\theta(x,y)$, $\theta \rightarrow \pi_{\theta}(x)$ and $\theta
\rightarrow \varphi_{x}(\theta)$, are continuous, and $\theta
\rightarrow f(s_0;\theta|x)$, as well as $\theta \rightarrow
f(s_1;\theta|x,s_0)$, are continuous for all $x \in {\mathcal  X}$ and
$s_0,s_1 \in {\mathbf R}^d$. Furthermore, for all $x \in {\mathcal  X}$
and $s_0,s_1 \in {\mathbf R}^d$, $f(s_0;\theta|x) \rightarrow 0$ and
$f(s_1;\theta|x,s_0) \rightarrow 0$, as $|\theta| \rightarrow \infty$.

C8. $E^{\theta_0}_x | \log (f(\xi_0;\theta_0|x)
f(\xi_1;\theta_0|x,\xi_0))| <  \infty$ for all $x \in {\mathcal  X}$.

C9. For each $\theta \in \Theta$, there is $\delta > 0$ such that, for
all $x \in {\mathcal  X}$,
\[
E^{\theta_0}_x  \biggl(\, \sup_{|\theta' - \theta| < \delta}  \bigl[\log
\bigl(f(\xi_0;\theta'|x)
 f(\xi_1;\theta'|x,\xi_0)\bigr) \bigr]^+  \biggr) < \infty,
\]
where $a^+ = \max\{a,0\}$. And there is a $b > 0$ such that, for all $x
\in {\mathcal  X}$,
\[
E_x^{\theta_0}  \biggl(\, \sup_{|\theta'| > b}  \bigl[\log
\bigl(f(\xi_0;\theta'|x)
 f(\xi_1;\theta'|x,\xi_0)\bigr) \bigr]^+  \biggr) < \infty.
\]
\begin{rem}\label{r2}
(a) Condition C1 is the $w$-uniform ergodicity condition for the
underlying Markov chain, which is considerably weaker than the
uniformly recurrent condition A1 of~\cite{JP99}, and that
of~\cite{DM01}. Furthermore, we impose conditions (\ref{5.1a}) and
(\ref{5.1b}) to guarantee that the induced Markovian iterated random
functions system satisfies Assumptions \hyperref[assK2]{K2} and
\hyperref[assK3]{K3} in Section~\ref{s4}.

(b) To have better understanding of these properties, we first consider
a simple state space model $X_n = \alpha X_{n-1} + \varepsilon_n,
\xi_n=X_n  + \eta_n$, where $|\alpha| < 1$, $\varepsilon_n$ and
$\eta_n$ are i.i.d. standard normal random variables, and they are
mutually independent. Since $\xi_n$ are independent for given $X_n$,
the weight function $w$ depends on $X_0$ only and we have $w(x)=|x|
+1$. Note that ${\mathcal  X} = {\mathbf R}$. Denote $b
=(1-|\alpha|^p)/(1-|\alpha|)$. Observe that
\begin{eqnarray*}
&&  \sup_{x \in {\mathbf R}}   \int_{- \infty}^{\infty}
 \frac{\exp\{ -(y - \alpha x)^2/2 \}}{\sqrt{2 \pi}}  \frac{\exp \{ -(s - y)^2/2 \} }{\sqrt{2 \pi}}\, dy
\\
&&\qquad =  \sup_{x \in {\mathbf R}} \frac{\sqrt{1/2}}{\sqrt{2 \pi}}
\exp \{ -(\alpha x - s)^2/4 \}
\\
&&\quad\qquad {}\times \int_{- \infty}^{\infty}
\frac{1}{\sqrt{2 \pi (1/2)}}  \exp\bigl\{ -\bigl(y - (\alpha x + s)/2\bigr)^2/2(1/2) \bigr\}\,  dy
\\
&&\qquad =  \frac{\sqrt{1/2}}{\sqrt{2 \pi}}\sup_{x \in {\mathbf R}}
\exp \{ -(\alpha x - s)^2/4 \} = \frac{1}{\sqrt{4\pi}}.
\end{eqnarray*}
A simple calculation leads to
\begin{eqnarray}\label{rem2.2}
\nonumber &&  \sup_{(x_0,s_0) \in {\mathbf R} \times {\mathbf R}}
E^\alpha_{(x_0,s_0)}  \biggl\{ \log  \biggl( g(s_0,\xi_1)^p
\frac{w(X_p,\xi_p)}{w(x_0,s_0)}   \biggr)  \biggr\}
\\
&&\qquad < \log  \sup_{x_0 \in {\mathbf R}} E^\alpha_{x_0} \biggl\{
\frac{|\alpha^p x_0 + \sum_{k=0}^{p-1} \alpha^k \varepsilon_{p-k}|+1}
{(4 \pi)^{p/2} (|x_0|+1)}  \biggr\}
\nonumber
\\[-8pt]
\\[-8pt]
\nonumber &&\qquad  \leq \log   \sup_{x_0 \in {\mathbf R}}  \biggl\{
\frac{|\alpha^p x_0| + E^\alpha_{x_0} |\sum_{k=0}^{p-1} \alpha^k
\varepsilon_{p-k}|+1} {(4 \pi)^{p/2}(|x_0|+1)}  \biggr\}
\\
\nonumber &&\qquad = \log \sup_{x_0 \in {\mathbf R}} \biggl\{
\frac{|\alpha^p x_0| + {2b}/{\sqrt{2\pi}}+1} {(4
\pi)^{p/2}(|x_0|+1)}  \biggr\} < 0.
\end{eqnarray}
This implies that (\ref{5.1a}) holds. By using the same argument, we
see (\ref{5.1b}) holds.

Next, we consider the case that $\varepsilon_n$ and $\eta_n$ are i.i.d.
double exponential$(1)$ random variables. Observe that
\begin{eqnarray*}
&& \sup_{x \in {\mathbf R}}  \int_{- \infty}^{\infty} \frac{\exp\{ -|y
- \alpha x| \}}{\sqrt{2}}  \frac{\exp \{ -|s - y| \} }{\sqrt{2}}\, dy
\\
&&\qquad = \frac{1}{4} \sup_{x \in {\mathbf R}} \bigl( (1 + |\alpha x
-s|) \exp \{ -|\alpha x - s| \} \bigr) = \frac{1}{4}.
\end{eqnarray*}
By making use of the same argument as in (\ref{rem2.2}), we see that
(\ref{5.1a}) and (\ref{5.1b}) hold. The extension to $\xi_n=
\beta_{X_n} \xi_{n-1} + \eta_n$, studied in Remark~\ref{r1}(b), is
straightforward and will not be repeated here. Other practical used
models of the Markov-switching model, \textup{ARMA} models,
(G)\mbox{ARCH} models and SV models will be given in
Section~\ref{s6}.

(c) Note that the mean contraction property ${\mathbf E} \log L_1 < 0$
is not satisfied in the above examples. Instead of applying
Theorem~\ref{T1} directly, we will explore the special structure of the likelihood
function in Lemma \ref{L4} below, such that $\{((X_n,\xi_n),M_n), n\geq
0\}$ is an irreducible, aperiodic and Harris recurrent Markov chain.
Hence, we can apply Theorem~\ref{T1} for the Markovian iterated functions
system on ${\mathbf M}$ induced from (\ref{2.4})--(\ref{2.7}).

(d) C2--C4 are standard smoothness conditions. C5 is the technical
condition for the existence of the Fisher information to be defined in
(\ref{5.6}) below. C8~and~C9 are integrability conditions that will be
used to prove strong consistency of the~MLE. Condition C6 is the
identifiability condition for state space models. That is, the family
of mixtures of $\{f(\xi_1;\theta|x,\xi_0)\dvtx  \theta \in \Theta\}$ is
identifiable. This condition will be used to prove strong consistency
of the MLE. Although it is difficult to check this condition in a
general state space model, in many models of interest the parameter
itself is identifiable only up to a permutation of states such as a
finite state hidden Markov model with normal distributions. A
sufficient condition for the identifiable issue can be found in
Theorem~\ref{T1} of~\cite{DM01}. See also the paper by It\^o, Amari and
Kobayashi \cite{IAK92} for necessary and sufficient conditions in the
case that the state space is finite and $\xi_i$ is a deterministic
function of~$X_i$.

(e) When the state space of the Markov chain $\{X_n, n \geq 0\}$ is
finite, and the observations $\xi_n$ are conditionally independent, this reduces to the
so-called hidden Markov model. It is easy to see that condition C1
implies (A1) by choosing $w(x)=1$, and conditions C2--C4 reduce to
(A2), (A3)~and~(A5) of~\cite{BRR98}. Conditions C6--C9 reduce to
conditions C1--C6 in \cite{L92}. We will discuss condition C5 in
Remark~\ref{r3} after Lemma~\ref{L5}.
\end{rem}

Let $\{(X_n,\xi_n), n \geq 0\}$ be the Markov chain defined in
(\ref{2.1}) and (\ref{2.2}). Recall from (\ref{2.8}) that the $\log$
likelihood can be written as
\begin{eqnarray}\label{5.2}
\nonumber l(\theta) &=& \log p_n(\xi_1,\ldots,\xi_n;\theta)
= \log \|{\mathbf P}_\theta(\xi_n) \circ \cdots \circ {\mathbf P}_\theta(\xi_1)
\circ {\mathbf P}_\theta(\xi_0) \pi\|
\\
&=& \log \frac{\|{\mathbf P}_\theta(\xi_n) \circ \cdots \circ {\mathbf
P}_\theta(\xi_1) \circ {\mathbf P}_\theta(\xi_0) \pi\|} {\|{\mathbf
P}_\theta(\xi_{n-1}) \circ \cdots \circ {\mathbf P}_\theta(\xi_1) \circ
{\mathbf P}_\theta(\xi_0) \pi\|}
\\
\nonumber &&{} + \cdots+\log \frac{\|{\mathbf
P}_\theta(\xi_1) \circ {\mathbf P}_\theta(\xi_0) \pi\|} {\|{\mathbf
P}_\theta(\xi_0) \pi\|}.
\end{eqnarray}
For each $n$, denote
\begin{equation}\label{5.3}
M_n:= {\mathbf P}_\theta(\xi_n) \circ \cdots \circ {\mathbf
P}_\theta(\xi_1) \circ {\mathbf P}_\theta(\xi_0)
\end{equation}
as the Markovian iterated random functions system on ${\mathbf M}$
induced from (\ref{2.4})--(\ref{2.7}). Then $\{((X_n,\xi_n),M_n), n\geq
0\}$ is a Markov chain on the state space $({\mathcal  X} \times
{\mathbf R}^d) \times {\mathbf M}$, with transition probability kernel
${\mathbf P}_\theta$ defined as in (\ref{3.2}). Let $\Pi_{\theta}$ be
the stationary distribution of $\{((X_n,\xi_n),M_n), n\geq 0\}$ defined
in Theorem~\ref{T1}(iii). Then the log-likelihood function $l(\theta)$
can be written as $S_n:=\sum_{k=1}^{n}g(M_{k-1},M_k)$ with
\begin{equation}\label{5.4}
g(M_{k-1},M_k):= \log \frac{\|{\mathbf P}_\theta(\xi_k) \circ \cdots
\circ {\mathbf P}_\theta(\xi_1) \circ {\mathbf P}_\theta(\xi_0) \pi\|}
{\|{\mathbf P}_\theta(\xi_{k-1}) \circ \cdots \circ {\mathbf
P}_\theta(\xi_1) \circ {\mathbf P}_\theta(\xi_0) \pi\|}.
\end{equation}

In order to apply Theorems \ref{T1}--\ref{T4}, we need to check that the Markovian
iterated random functions system satisfies Assumption K, and the
induced Markov chain is aperiodic, irreducible and Harris recurrent.
For this purpose, we need to define a suitable metric on the space
${\mathbf M}$, which has been defined in (\ref{2.4}). First, we add a
further condition on ${\mathbf M}$ to have
\[
{\mathbf M} = \biggl\{h|h\dvtx {\mathcal  X} \to {\mathbf R}^+ \mbox{
is }m\mbox{-measurable, }  \int h(x)m(dx) < \infty \mbox{ and }\sup_{x
\in {\mathcal X}} h(x)< \infty \biggr\}.
\]
For convenience of notation, we still use the notation $\mathbf M$, and
will use $h$ to represent an element in ${\mathbf M}$, which is
different from the notation $u$ used in Sections \ref{s3} and \ref{s4}.
We define the variation distance between any two elements $h_1,h_2$ in
${\mathbf M}$ by
\begin{equation}\label{5.5}
d(h_1,h_2) = \sup_{x \in \mathcal  X} |h_1(x)-h_2(x)|.
\end{equation}

Note that $({\mathbf M},d)$ is a complete metric space with Borel
$\sigma$-algebra ${\mathcal  B}({\mathbf M})$, but it is not separable.
Thus, Theorems \ref{T1}--\ref{T4} do not apply. However, rather than deal with the
measure-theoretic technicalities created by an inseparable space, we
can apply the results developed in Section~7 of~\cite{DF99} for
a direct argument of convergence. Therefore, Theorems \ref{T1}--\ref{T4} still hold
under the regularity conditions.

In order to describe our main results, we need the following lemmas
first. Their proofs are given in Section~\ref{s7}.

\begin{lemma}\label{L3}
Assume~\textup{C1--C5} hold or~\textup{C1, C6--C9} hold. Then for each
$\theta \in \Theta$ and $j=1,\ldots,n$, the random functions ${\mathbf
P}_{\theta}(\xi_0)$ and ${\mathbf P}_{\theta}(\xi_j)$, defined in
\textup{(\ref{2.5})}~and~\textup{(\ref{2.6})}, from $({\mathcal  X}
\times {\mathbf R}^d) \times {\mathbf M}$ to ${\mathbf M}$ are
Lipschitz continuous in the second argument, and the Markovian iterated
random functions system~\textup{(\ref{2.4})--(\ref{2.7})} satisfies
Assumption~\textup{\ref{assK}}. Furthermore, the function $g$ defined
in~\textup{(\ref{5.4})} belongs to ${\mathcal  L}^r(Q)$ for any $r >
0$.
\end{lemma}

For each $\theta \in \Theta$, recall that $\{((X_n,\xi_n),M_n),n \geq
0\}$ is a Markov chain induced by the Markovian iterated random
functions system~\textup{(2.4)--(2.7)} on the state space $({\mathcal
X} \times {\mathbf R}^d) \times {\mathbf M}.$

\begin{lemma}\label{L4}
Assume~\textup{C1--C5} hold or~\textup{C1, C6--C9} hold. Then for each
$\theta \in \Theta$, $\{((X_n,\xi_n),M_n),n \geq 0\}$ is an aperiodic,
$(m \times Q \times Q)$-irreducible and Harris recurrent Markov chain.
\end{lemma}

\begin{lemma}\label{L5}
Assume~\textup{C1--C5} hold. Then the Fisher information matrix
\begin{eqnarray}\label{5.6}
\nonumber {\mathbf I}(\theta) &=&  (I_{ij}(\theta)) \nonumber
\\
&=&  \biggl( {\mathbf E}_{\Pi}^{\theta}  \biggl[
  \biggl(\frac{\partial \log \|{\mathbf P}_{\theta}(\xi_1)\circ {\mathbf P}_{\theta}(\xi_0)\pi\|}{\partial \theta_i} \biggr)
\\
\nonumber &&\phantom{\biggl( {\mathbf E}_{\Pi}^{\theta}  \biggl[}
{}\times  \biggl(\frac{\partial \log \|{\mathbf P}_{\theta}(\xi_1)\circ
 {\mathbf P}_{\theta}(\xi_0)\pi \|}{\partial \theta_j}
 \biggr)  \biggr]  \biggr)
\end{eqnarray}
is positive definite for $\theta$ in a neighborhood
$N_\delta(\theta_0)$ of $\theta_0$. Recall that ${\mathbf
E}_{\Pi}^{\theta} := {\mathbf E}_{\Pi}$ is defined as the expectation
under ${\mathbf P}_{\Pi}$ in~\textup{(\ref{3.2})}.
\end{lemma}

\begin{rem}\label{r3}
Note that the Fisher information (\ref{5.6}) is defined as the expected
value under the stationary distribution $\Pi^{\theta}$ of the Markov
chain $\{((X_n,\xi_n), M_n), n \geq 0\}$. It is worth mentioning that only
$\xi_n$ appears in $M_n$, in which it reflects the nature of state
space models.

When the state space ${\mathcal  X}$ is finite, and the random
variables $\xi_n$ are\break conditionally independent for given~$X_n$, let $H:=
H(\xi_1,\xi_0,\xi_{-1},\ldots) =\break \sum_{m=-\infty}^1
H_m(\xi_1,\xi_0,\ldots),$ where
\begin{eqnarray*}
H_m(\xi_1,\xi_0,\ldots) &:=& E^{\theta^0} \biggl\{\frac{\partial  \log
 f(\xi_m;\theta|X_m)}{\partial \theta} \Big|\xi_1,\xi_0,\ldots\biggr\}
\\
&&{} - E^{\theta^0} \biggl\{\frac{\partial  \log
 f(\xi_m;\theta|X_m)}{\partial \theta} \Big|\xi_0,\xi_{-1},\ldots\biggr\}
\\
&&{} + E^{\theta^0} \biggl\{\frac{\partial  \log
 p_{\theta}(X_m, X_{m+1})} {\partial \theta} \Big|\xi_1,\xi_0,\ldots\biggr\}
\\
&&{} - E^{\theta^0} \biggl\{\frac{\partial  \log
 p_{\theta}(X_m X_{m+1})}{\partial \theta} \Big|\xi_0,\xi_{-1},\ldots\biggr\}.
\end{eqnarray*}
Under their Assumptions 1--4, Bickel and Ritov \cite{BR96} showed that
$H \in {\mathcal  L}^2(P^{\theta^0})$ and defined ${\mathbf
I_H}(\theta^{0}) := E^{\theta^0} \{H H^t \}.$ They also showed that
\[
\lim_{n \rightarrow \infty} \frac{1}{n}
 E^{\theta^0}  \biggl( \biggl(\frac{\partial \log\|T_n \pi
 \||_{\theta=\theta^0}} {\partial\theta}\biggr)
\biggl(\frac{\partial \log\|T_n \pi \||_{\theta=\theta^0}}{\partial
\theta}\biggr)^t  \biggr) = {\mathbf I_H}(\theta^0).
\]
In this paper we represent the log likelihood function of an additive
functional of the Markov chain $\{((X_n,\xi_n), M_n), n \geq 0\}$ in
(\ref{5.4}), and  then apply the strong law of large numbers for Markovian
iterated random functions given in Theorem~\ref{T1}(iv) to have, with
probability~1,
\[
 \lim_{n \rightarrow \infty} \frac{1}{n}
  \frac{\partial^2}{\partial\theta_i\, \partial \theta_j} \log\|{\mathbf P}_\theta(\xi_n)
  \circ \cdots \circ {\mathbf P}_\theta(\xi_1) \circ {\mathbf P}_\theta(\xi_0) \pi\| = -I_{ij}(\theta).
\]
Hence, under Assumptions 1--4 of~\cite{BR96}, ${\mathbf I}(\theta)$ is
well defined and is equal to ${\mathbf I_H}(\theta)$. The moment
condition in Assumption 4 of~\cite{BR96} can be relaxed to the
following: there exists a $\delta > 0$  with $ \rho_0 (\xi):=
\sup_{|\theta- \theta^{0}| < \delta} \max_{x,y
 \in {\mathcal  X}} \frac{f(\xi;\theta|x)}{f(\xi;\theta|y)},$ such that
$\sup_{x \in {\mathcal  X}} P^{\theta^0}\{\rho_0(\xi_1) = \infty |X_0
=x\} < 1$; see~\cite{BRR98}.
\end{rem}

\begin{lemma}\label{L6}
Assume~\textup{C1--C5} hold. Let $l'_j(\theta_0) = {\partial
l(\theta)}/ {\partial \theta_j}|_{\theta=\theta_0}$. Then,\break as $n
\rightarrow \infty,$
\begin{equation}\label{5.7}
 \frac1{\sqrt{n}}  ( l'_j(\theta_0)  )_{j=1,\ldots,q}
 \longrightarrow N\bigl(0, {\mathbf I}(\theta_0)\bigr)\qquad\mbox{in
 distribution.}
\end{equation}
\end{lemma}

\begin{theorem}\label{T5}
Assume~\textup{C1--C5} hold. Then there exists a sequence of solutions
$\hat{\theta}_n$ of~\textup{(\ref{5.1})} such that $\hat{\theta}_n \to
\theta_0$ in probability. Furthermore, $\sqrt{n}(\hat{\theta}_n -
\theta_0)$ is asymptotically normally distributed with mean zero and
variance--covariance matrix ${\mathbf I}^{-1}(\theta_0)$.
\end{theorem}

Since the proof of Theorem~\ref{T5} follows a standard argument, we
will not give it here.

\begin{corollary}\label{C2}
Under the assumptions of~Theorem~\textup{\ref{T5}}, if the likelihood equation
has a unique root for each $n$ and all $\xi_1,\ldots,\xi_n$, then there
is a consistent sequence of estimators $\hat{\theta}_n$ of the unknown
parameters $\theta_0$.
\end{corollary}

Next, we prove strong consistency of the MLE when the $\log$ likelihood
function is integrable. A crucial step is to give an appropriate
definition of the Kullback--Leibler information for state space models,
so that we can apply Theorem~\ref{T1} to have a standard argument of strong
consistency for the MLE. Here, we define the Kullback--Leibler
information as
\begin{eqnarray}\label{5.8}
K(\theta_0,\theta) &=& {\mathbf E}_{\Pi}^{\theta_0}  \biggl(
\log\frac{\| {\mathbf P}_{\theta_0}(\xi_1) \circ {\mathbf
P}_{\theta_0}(\xi_0) \pi_{\theta_0}\|} {\|{\mathbf P}_{\theta}(\xi_1)
\circ {\mathbf P}_{\theta}(\xi_0) \pi_{\theta}\|}  \biggr) \nonumber
\\[-8pt]
\\[-8pt]
\nonumber &\hspace*{-3pt} :=& \int  \log\frac{\| {\mathbf P}_{\theta_0}(\xi_1) \circ
{\mathbf P}_{\theta_0}(\xi_0) \pi_{\theta_0}\|} {\|{\mathbf
P}_{\theta}(\xi_1) \circ {\mathbf P}_{\theta}(\xi_0) \pi_{\theta} \|}
\Pi\bigl(d(x,\xi)\times d \pi_{\theta_0}\bigr).
\end{eqnarray}

\begin{theorem}\label{T6}
Assume that~\textup{C1, C6--C9} hold and let $\hat{\theta}_n$ be the
MLE based on n observations $\xi_0,\xi_1,\ldots,\xi_n$. Then
$\hat{\theta}_n \longrightarrow \theta_0$ $P^{\theta_0}$-a.s. as
$n\rightarrow\infty$.
\end{theorem}

Since the proof of Theorem~\ref{T6} follows a standard argument, we
will not give it here.

To derive the  Edgeworth expansion for the MLE, we need to define the
following notation and assumptions first. For nonnegative integral
vectors $\nu=(\nu^{(1)},\ldots,\nu^{(q)})$, write
$|\nu|=\nu^{(1)}+\cdots+\nu^{(q)}$, $\nu!=\nu^{(1)}! \cdots
\nu^{(q)}!$, and let $D^{\nu}=(D_1)^{\nu^{(1)}} \cdots
(D_q)^{\nu^{(q)}}$ denote the $\nu$th derivative with respect to
$\theta$. Suppose assumptions C2, C3, C4 and C5 are strengthened so
that there exists~$r \geq 3$, as follows.

C2$'$.  The true parameter $\theta_0$ is an interior point of $\Theta$.
For all $x\in {\mathcal  X}$, $s_0,s_1 \in {\mathbf R}^d$, $\theta \in
\Theta \subset {\mathbf R}^q$, the partial derivatives
\[
D^1 f(s_0;\theta|x),\qquad  D^2 f(s_0;\theta|x),\ldots, D^r
f(s_0;\theta|x),
\]
as well as the partial derivatives
\[
D^1 f(s_1;\theta|x,s_0),\qquad D^2 f(s_1;\theta|x,s_0),\ldots,  D^r
f(s_1;\theta|x,s_0),
\]
and for all $x,y \in {\mathcal X}$, $\theta \rightarrow p_\theta(x,y)$
and $\theta \rightarrow \pi_{\theta}(x)$ have $r-1$ continuous
derivatives in some neighborhood $N_\delta(\theta_{0}):=\{\theta\dvtx
|\theta-\theta_0| < \delta\}$ of $\theta_0$.

C3$'$.
\[
\int_{\mathcal  X} \sup_{\theta \in N_\delta(\theta_0)} |
D^1 \pi_\theta(x)| m(dx) < \infty,\ldots,  \int_{\mathcal  X}
\sup_{\theta \in N_\delta(\theta_0)} |D^{r-1}  \pi_{\theta}(x)|
m(dx)<\infty,
\]
 and for all $x \in {\mathcal  X}$,
\[
\int_{\mathcal  X} \sup_{\theta \in N_\delta(\theta_0)} |D^1
p_{\theta}(x,y)| m(dy) < \infty,\ldots, \int_{\mathcal  X} \sup_{\theta
\in N_\delta(\theta_0)} | D^{r-1}  p_{\theta}(x,y)| m(dy) < \infty.
\]

C4$'$. For all $x \in {\mathcal  X}$, $s_0 \in {\mathbf R}^d$ and
$\theta \in \Theta$,
\[
E_x^{\theta}\big| D^{\nu} f(\xi_0;\theta|x)\big|^r < \infty,\qquad
E_{(x,s_0)}^{\theta}\big| D^{\nu}f(\xi_1;\theta|x,s_0)\big|^r < \infty,
\]
for $1 \leq |\nu| \leq r$, and
\begin{eqnarray*}
E_x^{\theta} \biggl(\, \sup_{\theta \in N_\delta(\theta_0)} \big| D^{\nu}
f(\xi_0;\theta|x)\big|^r \biggr) &<& \infty,
\\
E_{(x,s_0)}^{\theta}
\biggl(\, \sup_{\theta \in N_\delta(\theta_0)} \big|
D^{\nu}f(\xi_1;\theta|x,s_0)\big|^r \biggr) &<& \infty,
\end{eqnarray*}
for $|\nu| = r+1$.

C5$'$.
\[
\sup_{x \in {\mathcal  X}} E_{x}^{\theta_0}  \biggl(\, \sup_{|\theta-
\theta_0| < \delta} \sup_{y,z \in {\mathcal  X}}
\frac{f(\xi_0;\theta|y)f(\xi_1;\theta|y, \xi_0)}
{f(\xi_0;\theta|z)f(\xi_1;\theta|z,\xi_0)}  \biggr)^r  < \infty.
\]

We will assume conditions (\ref{4.20}) and (\ref{4.21}) hold for
$Z_j^{(\nu)}:= D^{\nu} \log p_1(\xi_0,\break \xi_1;\theta_0)$, $1 \leq |\nu|
\leq r$. Let $Z_j:=\{Z_j^{(\nu)}\dvtx  1 \leq |\nu| \leq r\}$ be
$p$-dimensional random vectors for $j \geq 1$, where $p$ is the number
of all distinct multi-indices $\nu$, \mbox{$1 \leq |\nu| \leq r$}. In the
following, denote $\bar{Z}=(1/n)\sum_{k=1}^n Z_k$.

Use a standard argument involving the sign change of a continuous
function, or a fixed point theorem in the multi-parameter case
(cf.~\cite{BG78}), to prove that the likelihood equation has a solution
which converges in probability to $\theta_0$. Note that the following
notation is interpreted in the multi-dimensional sense. Applying the
moderate deviation result on $\bar{Z}$ in Theorem~\ref{T4}, it is possible to
ensure that, with $P^{\theta_0}$-probability $1-o(n^{-1})$,
$\hat{\theta}_n$ satisfies the likelihood equation and lies on
$(\theta_0\pm\log n/\sqrt{n}\,)$. It is this solution we take as our
$\hat{\theta}_n$. If the likelihood equation has multiple roots, assume
we have a consistent estimator $T_n$ such that $T_n$ lies in
$(\theta_0\pm\log n/\sqrt{n}\,)$ with $P^{\theta_0}$-probability
$1-o(n^{-1})$. In this case, we may take the solution nearest to $T_n$.
By the preceding reasoning, this solution, which is identifiable from
the sample, will lie in $(\theta_0\pm\log n/\sqrt{n}\,)$ with
$P^{\theta_0}$-probability $1-o(n^{-1})$.

Clearly, with $\hat{\theta}_n$ as above, with probability
$1-o(n^{-1})$,
\begin{equation}\label{5.9}
\hspace*{12pt}
0= \bar{Z}^{(e_s)} + \sum_{|\nu|=1}^{r-1} \frac{1}{\nu!}
\bar{Z}^{(e_s+\nu)}
(\hat{\theta}_n-\theta_0)^{\nu}+R_{n,s}(\hat{\theta}_n),\qquad 1 \leq s
\leq q,
\end{equation}
where $e_s$ has $1$ as the $s$th coordinate and zeros otherwise.

We rewrite equation (\ref{5.9}) as
\begin{equation}\label{5.10}
0=A(\bar{Z},\hat{\theta}_n)+R_n.
\end{equation}
Note $0=A(\gamma(\theta_0),\theta_0)$ and
$ \frac{\partial A}{\partial \theta} |_{\gamma(\theta_0),\theta_0}=
-(\mbox{Fisher information})\neq 0.$

Hence, by the implicit function theorem, there are a neighborhood $N$ of
$\gamma$ and $q$ uniquely defined real-valued infinitely differentiable
functions ${\mathbf g}_i$ ($1 \leq i \leq q$) on $N$ such that $\theta =
{\mathbf g}(z) = ({\mathbf g}_1(z),\ldots,{\mathbf g}_q(z))$ satisfies
(\ref{5.10}). This implies, with probability $1-o(n^{-1})$,
$|\hat{\theta}_n-\theta_0|\leq K(\log n/\sqrt{n}\,)^4.$

To derive the asymptotic expansion of
$P^{\theta_0}\{\sqrt{n}(\hat{\theta}_n- \theta_0) \in B \}$, note that
$\hat{\theta}_n= {\mathbf g}(n^{-1}Z_n)$, where ${\mathbf g}\dvtx
{\mathbf R}^p \rightarrow {\mathbf R}^q$ is sufficiently smooth in some
neighborhood of $\gamma$. For the case of i.i.d. $\xi_n$, Bhattacharya
and Ghosh \cite{BG78} made use of the Edgeworth expansion of the
distribution of $(S_n-n\gamma)/\sqrt{n}$ to derive an Edgeworth
expansion of the distribution of $\sqrt{n}\{{\mathbf
g}(n^{-1}S_n)-{\mathbf g}(\gamma)\}$. Making use of Theorem~\ref{T4} and a
straightforward extension of their argument, we can generalize their
result to have the following theorem.

\begin{theorem}\label{T7}
Assume~\textup{C1, C2$'$--C5$'$} hold for some $r \geq 3$. Assume
\textup{(\ref{4.20})} and~\textup{(\ref{4.21})} hold. Let $J_{\mathbf
g}=(D_j {\mathbf g}_i(\gamma))_{1\leq i \leq q,1\leq j \leq p}$ be the
$q \times p$ Jacobian matrix and let $V({\mathbf g}) = J_{\mathbf g}V
J_{\mathbf g}'$. Then there exists a sequence of solutions
$\hat{\theta}_n$ of~\textup{(\ref{5.1})}, and there exist polynomials
$p_j$ in $q$ variables $(1 \leq j \leq r-2)$ such that
\begin{eqnarray*}
&& \sup_{B\in {\mathcal  B}_{a,c}}\Bigg|{\mathbf
P}^{\theta_0}_\nu\bigl\{\sqrt{n}(\hat{\theta}_n - \theta_0) \in B
\bigr\} - \int_B \Biggl\{\phi_{V({\mathbf g})}(y)+
\sum_{j=1}^{r-2}n^{-j/2}\phi_{j,V,{\mathbf g}}(y)\Biggr\}\,dy\Bigg|
\\
&&\qquad =
o\bigl(n^{-(r-2)/2}\bigr),
\end{eqnarray*}
where $\phi_{j,V,{\mathbf g}} = \tilde{\pi}_{j,{\mathbf g}}(-D)\phi_V $
and $\tilde{\pi}_{j,{\mathbf g}}(y)$ is a polynomial in $y (\in
{\mathbf R}^p)$ whose coefficients are smooth functions of the partial
derivatives of $\lambda(\alpha)$ at $\alpha =0$ up to order $j+2$, and
those of $\nu_\alpha {\mathbf Q}_\alpha h_1$ at $\alpha = 0$ up to
order $j$ together with those of ${\mathbf g}$ at $\mu$ up to order~$j
+1$.
\end{theorem}

The application of Theorem~\ref{T7} to third-order efficiency for the MLE and
third-order efficient approximate solution of the likelihood equation
follows directly from~\cite{G94}.

\section{Examples}\label{s6}
From a theoretical point of view, Theorems \ref{T5}--\ref{T7} are adequate for state
space model estimation problems in providing assurance of the
existence of efficient estimators, characterizing them as solutions of
likelihood equations and prescribing their asymptotic behavior. In
practice, however, one must still contend with certain statistical and
numerical difficulties, such as implementation of the maximum
likelihood estimator. In this section we apply our results to study
some examples which include Markov switching models \textup{ARMA}
models, (G)\mbox{ARCH} models and SV models. For simplicity, in
these examples we consider only specific structure of normal error
assumption in most cases. Although strong consistency and asymptotic
normality of the MLE in \textup{ARMA} and $\mbox{GARCH}(p,q)$ have
been known in the literature, we provide alternative proofs in the
framework of state space models. Furthermore, we can apply Theorem~\ref{T7} to
have Edgeworth expansion for the MLE. To the best of our knowledge, the
asymptotic normality of the MLE in the $\mbox{AR}(1)/\mbox{ARCH}(1)$ model,
considered in Section~\ref{s6.3}, seems to be new. The results of
asymptotic properties for the MLE in stochastic volatility models not
only provide theoretical justification, but also give some insight into the
structure of the likelihood function, which can be used for
further study.

\subsection{Markov switching models}

We start with a simple real-valued fourth-order autoregression around
one of two constants, $\mu_1$ or $\mu_2$:
\begin{equation}\label{mg1}
 \xi_n - \mu_{X_n}
= \sum_{k=1}^4 \varphi_k\bigl(\xi_{n-k} - \mu_{X_{n-k}}\bigr) +
\varepsilon_n,
\end{equation}
where $\varepsilon_n \sim N(0,\sigma^2)$, and $\{X_n, n \geq 0\}$ is a
two-state Markov chain. This model was studied by Hamilton \cite{H89}
in order to analyze the behavior of U.S. real GNP. To apply our theory
in the form of (\ref{mg1}), we consider a simple case of order $1$ in
(\ref{mg1}). In this case, the likelihood function for given $X_n =
x_n$, $n\geq 0$, is
\begin{equation}\label{mg2}
\hspace*{12pt}
f(\xi_n|x_n;\theta) = \frac{1}{\sqrt{2\pi}\sigma} \exp  \bigl( -\bigl[
\bigl(\xi_n - \mu_{x_n}\bigr) -  \varphi_1\bigl(\xi_{n-1} -
\mu_{x_{n-1}}\bigr)\bigr]^2/(2\sigma^2)  \bigr).
\end{equation}
Denote by $[p_{xy}]_{x,y=1,2}$ the transition probability of the
underlying Markov chain $\{X_n, n \geq 0\}$ and let $\theta =(p_{11},
p_{21},\varphi_1,\mu_1,\mu_2,\sigma^2)$ be the unknown parameter.
Assume that $| \varphi_1 | <1$, and that there exists a constant $c>0$
such that $\sigma^2 >c$. Moreover, we assume that $\mu_1 \neq \mu_2$
such that the identifiability condition C6 holds. Since the state space
of $X_n$ is finite, we consider $0< p_{xy} < 1$ for all $x,y =1,2$, and
let $w(x)= |x| +1 $ such that the condition C1 holds. Under the normal
distribution assumption, it is easy to see that conditions C2--C4 and
C7--C9 are satisfied in this model. To check that C5 holds note that
condition C5 reduces to
\begin{equation}\label{mg3}
\sup_{x \in {\mathcal  X} } E_{x}^{\theta^0}  \biggl(  \biggl[
\sup_{|\theta- \theta^{0}| < \delta} \max_{y, z \in {\mathcal  X}}
\frac{f(\xi_0;\theta|y)f(\xi_1;\theta|y,\xi_0)}{f(\xi_0;\theta|z)f(\xi_1;\theta|z,\xi_0)}
 \biggr]^2   \biggr) < \infty.
\end{equation}
Since the maximum over $x$, $y$ and $z$ is applied to a finite set
${\mathcal  X}$, and $f$ defined in (\ref{mg1}) is a normal density, it
is easy to check that (\ref{mg3}) is satisfied.

When $\xi_n=X_n$ as in (\ref{mg1}), that is, $\mu_1=\mu_2=\mu$ are
given, this reduces to the classical autoregressive model with unknown
parameters $\theta=(\varphi_1,\ldots,\varphi_4,\sigma^2)$. The Fisher
information matrix is then given by
\begin{equation}\label{mg4}
 {\mathbf I}(\theta)=   \pmatrix{
 \sigma^{-2} \Gamma & 0  \cr\noalign{}
 0  &  2(\sigma^4)^{-1}},
\end{equation}
where $\Gamma = (\gamma_{i-j})_{4\times 4}$ for $1 \leq i, j \leq 4$
with $\gamma_k = E X_n X_{n + k}$. A simple calculation shows that
(\ref{5.6}) reduces to (\ref{mg4}) in this case. When $\varphi_k=0$ as
in (\ref{mg1}), this is the hidden Markov model with normal mixture
distributions considered in Example~1 of~\cite{BRR98}.

\subsection{ARMA models}

We start with a univariate Gaussian causal $\mbox{ARMA}(p,q)$ model
which can be written as a state space model by defining
$r=\max\{p,q+1\}$,
\begin{eqnarray}\label{6.1}
\xi_n - \mu &=& \alpha_1(\xi_{n-1} - \mu) + \alpha_2(\xi_{n-2} - \mu)+
\cdots+ \alpha_r(\xi_{n-r} - \mu) \nonumber
\\[-8pt]
\\[-8pt]
\nonumber &&{}  +  \varepsilon_n + \beta_1 \varepsilon_{n-1} + \beta_2
\varepsilon_{n-2}+\cdots+
  \beta_{r-1} \varepsilon_{n-r+1},
\end{eqnarray}
where $\alpha_j = 0$ for $j > p$ and $\beta_j = 0$ for $j > q$.
Furthermore, we assume $\varepsilon_n$ are i.i.d. random variables with
distribution $N(0,\sigma^2)$. Asymptotic properties of the MLE in
the \textup{ARMA} model can be found in \cite{H73} and
\cite{YB01}. A general treatment of the MLE in the  Gaussian
$\mbox{ARMAX}$ model can be found in Chapter~7 of~\cite{C88}.

By using the same idea as that in \cite{H94}, we consider the
following state space representation of (\ref{6.1}):
\begin{equation}\label{6.2}
X_{n+1} =  \left[\matrix{ \alpha_1 & \alpha_2 & \cdots & \alpha_{r-1} &
\alpha_r \cr\noalign{}
 1 & 0 &\cdots& 0 & 0 \cr\noalign{}
 0 & 1 &\cdots& 0 & 0 \cr\noalign{}
\vdots & \vdots& \ddots & \vdots & \vdots\cr\noalign{} 0 &
0 &\cdots& 1 & 0}\right] X_n + \left[\matrix{\varepsilon_{n+1}
\cr\noalign{}  0  \cr\noalign{} 0
\cr\noalign{}  \vdots \cr\noalign{}
0}\right]
\end{equation}
and
\begin{equation}\label{6.3}
\xi_{n} = \mu +  [ \begin{array}{ccccc} 1& \beta_1 &\beta_2  \cdots
\beta_{r-1} \end{array}  ] X_n.
\end{equation}

Assume that the roots of $1-\alpha_1z - \alpha_2z^2 - \cdots - \alpha_p
z^p = 0$ lie outside the unit circle. It is easy to see that $\{X_n, n
\geq 0 \}$ forms a $w$-uniformly ergodic Markov chain with
$w(x)=\|x\|^2$ (cf. Theorem~16.5.1 in \cite{MT93}). And $\xi_n$ are
conditionally independent given $\{X_n, n \geq 0\}$. Since the
verification of the weighted mean contraction property and the weighted
moment assumption is the same as those in Remark~\ref{r2}(b), it will not be
repeated here. This implies that condition C1 holds. The assumption
$\varepsilon_n \sim N(0,\sigma^2)$ also implies that conditions C2--C5,
C2$'$--C5$'$ and C7--C9 are satisfied in model (\ref{6.1}). Since the
verification is straightforward, we do not report it here. Suppose the
conditional distribution of $\xi_n$ given $X_0,\ldots,X_n$ is of the
form $F_{X_{n-1},X_n}$ from (\ref{6.3}). The Cram{\'e}r conditions
(\ref{4.20})~and~(\ref{4.21}) hold for $Z_j^{(\nu)}:= D^{\nu} \log
p_1(\xi_0,\xi_1;\theta_0)$, since the conditional density of $\xi_n$
given $\{x_n, n \geq 0\}$ is $N(0,\sigma^2)$ and
\begin{equation}\label{6.4}
\limsup_{|\theta| \to 0} \bigg|\int_{-\infty}^{\infty}
\int_{-\infty}^{\infty} \biggl\{\int_{-\infty}^{\infty} e^{i\theta \xi}
\,d F_{x,\alpha x+z} (\xi)\biggr\} \varphi(z)\,dz \,\pi(dx)\bigg| < 1,
\end{equation}
where $\varphi(\cdot)$ is the normal density function of
$\varepsilon_1$, and $\pi$ is the stationary distribution of $\{X_n\}$.
The identification issue in C6 can be found in Chapter~9 of~\cite{C88}
or Chapter~13 of~\cite{H94}.

\subsection{$(G)\mbox{ARCH}$ models}\label{s6.3}

In this subsection we study two specific $(\mbox{G})\mbox{ARCH}$ models. To
start with, we consider the $\mbox{AR}(1)/\mbox{ARCH}(1)$ model
\begin{equation}\label{6.5}
X_n = \beta_0 + \beta_1 X_{n-1} + \sqrt{\alpha_0 + \alpha_1 X_{n-1}^2}
\varepsilon_n,
\end{equation}
where $\alpha_i,\beta_i$ are unknown parameters for $i=0,1$ with
$\alpha_0 > 0, 0 < \alpha_1 < 1$, $3\alpha_1^2 < 1$ and $0< \beta_1 <
1$. Here $\varepsilon_n$ are i.i.d. random variables with the standard
normal distribution. Note that in (\ref{6.5}) ${\mathbf X}=(X_n)$ is
defined as the autoregressive scheme $\mbox{AR}(1)$ with $\mbox{ARCH}(1)$
noise $(\sqrt{\alpha_0 + \alpha_1 X_{n-1}^2} \varepsilon_n)_{n \geq
1}$. When $\beta_0=\beta_1=0$, this is the classical $\mbox{ARCH}(1)$
model first considered by Engle \cite{E82}.

Model (\ref{6.5}) is conditionally Gaussian, and therefore the likelihood
function of the parameter $\theta=(\alpha_0,\alpha_1,\beta_0,\beta_1)$
for given observations ${\mathbf x}=(x_0=0,x_1,\ldots,x_n)$ from~(\ref{6.5})~is
\begin{eqnarray}\label{6.6}
l({\mathbf x};\theta) &=&  (2 \pi)^{-n/2} \prod_{k=1}^n (\alpha_0 +
\alpha_1 x_{k-1}^2)^{-1/2}
\nonumber
\\[-8pt]
\\[-8pt]
\nonumber
&&\phantom{(2 \pi)^{-n/2} \prod_{k=1}^n }
{}\times \exp  \Biggl\{ - \frac{1}{2} \sum_{k=1}^n
\frac{(x_k - \beta_0 - \beta_1 x_{k-1})^2} {\alpha_0 + \alpha_1
x_{k-1}^2} \Biggr\}.
\end{eqnarray}

Assume $\beta_0=0$ and $\alpha_0,\alpha_1$ are given. The maximum
likelihood estimator $\hat{\beta}_1$ of $\beta_1$ is the root of the
equation $\partial l({\mathbf x};\theta)/\partial \beta_1 =  0.$
In view of (\ref{6.5}) and (\ref{6.6}), we obtain
\begin{eqnarray}\label{6.7}
 \hat{\beta}_1 &=&  \frac{{\sum_{k=1}^n (x_k - \beta_0) x_{k-1}/(\alpha_0 + \alpha_1 x_{k-1}^2)}}
{{\sum_{k=1}^n x_{k-1}^2/ (\alpha_0 + \alpha_1 x_{k-1}^2)}} \nonumber
\\[-8pt]
\\[-8pt]
\nonumber &=& \beta_1 + \frac{{\sum_{k=1}^n
x_{k-1}\varepsilon_{k}/\sqrt{\alpha_0 + \alpha_1 x_{k-1}^2}}}{
{\sum_{k=1}^n x_{k-1}^2/(\alpha_0 + \alpha_1 x_{k-1}^2)}}.
\end{eqnarray}

Meyn and Tweedie \cite{MT93}, pages 380 and 383, establish $w$-uniform
ergodicity [with $w(x)=|x|+1$] of the $\mbox{AR}(1)$ model $X_n = \beta_0 +
\beta_1 X_{n-1} + \varepsilon_n$ by proving that a \textit{drift
condition} is satisfied, where $|\beta_1| < 1$ and the $\varepsilon_n$
are i.i.d. random variables, with $E|\varepsilon_n| < \infty$, whose
common density function $q$ with respect to Lebesgue measure is
positive everywhere. The strongly nonlattice condition holds as that in
model (\ref{6.1}). By using an argument similar to Theorem~1
of~\cite{L96}, we have the asymptotic identifiability of the likelihood
function (\ref{6.6}). Letting $\xi_n=X_n$, and using an argument
similar to that in Remark~\ref{r2}(b), condition C1 holds. The verification of
conditions C2--C9 and C2$'$--C5$'$ is straightforward and tedious, and
is thus omitted. By Theorems \ref{T5}--\ref{T7}, we have the strong consistency,
asymptotic normality and Edgeworth expansion of the MLE
$\hat{\beta}_1$. The asymptotic properties of the MLE of $\beta_0$,
$\alpha_0$ and $\alpha_1$ can be verified in a similar way.

Next, we consider the $\mbox{GARCH}(p,q)$ model of (\ref{6.8}) in
Example~\ref{e1}. It is known that the necessary and sufficient condition
for (\ref{6.8}) defining a unique strictly stationary process $\{Y_n, n
\geq 0\}$ with $E Y_n^2 < \infty$ is
\begin{equation}\label{6.9}
\sum_{i=1}^p \alpha_i  + \sum_{j=1}^q \beta_j < 1.
\end{equation}
We assume (\ref{6.9}) holds.

Similar to the estimation for \textup{ARMA} models, the most
frequently used \mbox{estimators} for $\mbox{GARCH}$ models are those
derived from a (conditional) Gaussian likelihood function
(cf.~\cite{FY03}). Without the normal assumption of $\varepsilon_n$ in
(\ref{6.8}), and imposing the moment condition $E(\varepsilon_1^4) <
\infty$, Hall and Yao \cite{HY03} established the asymptotic normality
of the conditional maximum likelihood estimator in
$\mbox{GARCH}(p,q)$. They also established asymptotic results when
the case of the error distribution is heavy-tailed. Earlier in the
literature, when $p=q=1$, Lee and Hansen \cite{LH94} and
\mbox{Lumsdaine} \cite{L96} proved, under some regularity conditions, the consistency
and asymptotic normality for the quasi-maximum likelihood estimator in
the  $\mbox{GARCH}(1,1)$ model.\looseness=1

By using the state space representation (\ref{6.10}) and (\ref{6.11}),
it is known (cf. Theorem~3.2 of~\cite{BR92}) that the Markov chain
$\{X_n, n \geq 0 \}$ defined in (\ref{6.11}) is stationary if and only
if the top Lyapunov exponent $\gamma$ of $A_n$ is strictly negative. It
is easy to see that $\{X_n, n \geq 0\}$ is an aperiodic, irreducible
and $w$-uniformly [with $w(x)=\|x\|^2$] ergodic Markov chain.
Furthermore, we assume $\varepsilon_n$ are i.i.d. random variables with
distribution $N(0,\sigma^2)$. An argument similar to that in
Remark~\ref{r2}(b) leads to condition C1 holding. The normal error assumption
also implies that conditions C2--C5, \mbox{C2$'$--C5$'$} and \mbox{C7--C9} are
satisfied in model (\ref{6.11}). When $p=q=1$, Theorem~1 of~\cite{L96}
proves the asymptotic identifiability of the likelihood function.

\subsection{Stochastic volatility models}

Consider the stochastic volatility model (\ref{6.12})--(\ref{6.18}).
To check that condition C1 holds, 
we note that $w(x)=|x|+1$ in the $\mbox{AR}(1)$ model $X_{n} = \alpha X_{n-1}
+ \eta_n$ by proving that a \textit{drift condition} is satisfied,
where $|\alpha| < 1$ and the $\eta_n$ are i.i.d. random variables, with
$E|\eta_1| < \infty$, whose common density function $q$ with respect to
Lebesgue measure is positive everywhere. Since $\varepsilon_n \sim
N(0,1)$, $\zeta_n = \log \varepsilon^2$, $\eta_n \sim
N(0,\sigma_{\eta}^2)$, and $\zeta_n$ and $\eta_n$ are mutually
independent, an argument similar to that in Remark~\ref{r2}(b) leads to the result that the
rest of condition C1 holds. Conditions C2--C5, C2$'$--C5$'$ and
\mbox{C7--C9} are also satisfied in model (\ref{6.15})~and~(\ref{6.16}) (cf. pages
22--23 of~\cite{S96}). Denote $\xi_n:=\log Y_n^2$. Note that the
conditional density of $X_n$ exists, and this implies that the
conditional distribution of $\xi_n$ given $X_0,\ldots,X_n$ is of the
form $F_{X_{n-1},X_n}$ such that
\begin{equation}\label{6.19}
\limsup_{|t| \to 0} \bigg| \int_{\mathcal  X} \int_{-\infty}^{\infty}
\biggl\{\int_{-\infty}^{\infty} e^{i t s} \,d F_{x,\alpha x+z}
(s)\biggr\} \varphi(z)\,dz \,\pi(dx)\bigg| < 1,
\end{equation}
where $\varphi(\cdot)$ is the normal density function of $\zeta_1$ and
$\pi$ is the stationary distribution of $\{X_n\}$. Let
$S_n=\sum_{i=1}^n \xi_i$, $S_0=0$. Then $\{(X_n,S_n), n \geq 0\}$ is
strongly nonlattice. To check the identification condition C6, the
reader is referred to  Chapter~13 of~\cite{H94} and Section~2.4.3
of~\cite{CHR96}.

Next, we assume that $\varepsilon_n \sim N(0,1)$, $\zeta_n=\log
\varepsilon^2_n$ and $\eta_n$ is a sequence of i.i.d. double
exponential$(1)$ random variables. Furthermore, we assume $\zeta_n$ and
$\eta_n$ are mutually independent. By using an argument similar to that
in Remark~\ref{r2}(b), condition C1 holds. Simple calculations also lead
conditions C2--C5, C2$'$--C5$'$ and \mbox{C7--C9} to hold in this case. Under
the assumption that the conditional distribution of $\xi_n$ given
$X_0,\ldots,X_n$ is of the form $F_{X_{n-1},X_n}$ such
that~(\ref{6.19}) holds, $\{(X_n,S_n), n \geq 0\}$ is strongly
nonlattice.

Without the normal assumption, quasi-maximum likelihood (QML)
estimators of the parameters are obtained by treating $\zeta_n$ and
$\eta_n$ as though they were normal and maximizing the prediction error
decomposition form of the likelihood obtained via the Kalman filter or
implied volatility. That is, we assume that $\zeta_n$ is a sequence of
independent and identically distributed
$N(0,\sigma_{\zeta}^2)$ random variables. For given observations ${\mathbf y}=(\log
y_1^2,\ldots,\log y_n^2)$ from (\ref{6.15}) and (\ref{6.16}), the
likelihood function of the parameter $\theta=(\alpha,\sigma_\eta^2,
\sigma_\zeta^2)$ is
\begin{eqnarray}
\nonumber \hspace*{2mm}  l({\mathbf y};\theta) &=&  \int_{x_0 \in {\mathcal  X}} \cdots \int_{x_n
\in {\mathcal  X}} \pi(x_0) (2 \pi \sigma_\zeta^2)^{-n/2}
\\
&&\phantom{\int_{x_0 \in {\mathcal  X}} \cdots \int}
{}\times \prod_{k=1}^n p(x_{k-1},x_k)
\\
\nonumber &&\phantom{\int_{x_0 \in {\mathcal  X}} \cdots \int \times \prod_{k=1}^n}
{}\times \exp  \Biggl\{ - \frac{1}{2} \sum_{k=1}^n \frac{(\log y_k^2 -
\omega - x_{k})^2} {\sigma_\zeta^2} \Biggr\}\,dx_0\,dx_1 \cdots dx_n,
\end{eqnarray}
where $p(x_{k-1},x_k)$ is defined in (\ref{6.17}). By using the results
of~\cite{D79}, Harvey, Ruiz and Shephard \cite{HRS94} showed that the quasi-maximum likelihood estimators
are asymptotically normal under some regularity conditions. Further
study of the MLE in stochastic volatility models will be published in a
separate paper.

\section{\texorpdfstring{Proofs of Lemmas \protect\ref{L3}--\protect\ref{L6}}{Proofs of
Lemmas 3--6}}\label{s7}
 For convenience of notation, denote
$\{Z_n, n \geq 0\}:= \{((X_n,\xi_n),M_n),n \geq 0\}$ as the Markov
chain induced by the Markovian iterated random functions system
\textup{(\ref{2.4})--(\ref{2.7})} on the state space $({\mathcal X}
\times {\mathbf R}^d) \times {\mathbf M}.$ In the proof of Lemma~\ref{L3}, we
omit $\theta$ in ${\mathbf P}_\theta (\cdot)$ for simplicity.

\begin{pf*}{Proof of Lemma~\ref{L3}}
We consider only the cases of ${\mathbf P}(\xi_1)$, since the cases of
${\mathbf P}(\xi_0)$ and ${\mathbf P}(\xi_j)$, for $j=2,\ldots,n$, are a
straightforward consequence. For any two elements $h_1,h_2 \in {\mathbf
M}$, and two fixed elements $s_0, s_1 \in {\mathbf R}^d$,
by~(\ref{5.5}) we have
\begin{eqnarray*}
&&  d\bigl({\mathbf P}(s_1)h_1,{\mathbf P}(s_1)h_2\bigr)
\\
&&\qquad = \sup_{x_0 \in {\mathcal  X}} \bigg|\int p_\theta(x_0,x_1)
f(s_1;\theta|x_1,s_0) h_1(x_1) m(d x_1)
\\
&&\phantom{\qquad = \sup_{x_0 \in {\mathcal  X}} \bigg|}
{}- \int
p_\theta(x_0,x_1) f(s_1;\theta|x_1,s_0) h_2(x_1) m(d x_1)\bigg|
\\
&&\qquad \leq d(h_1,h_2) \sup_{x_0 \in {\mathcal  X}} \int
p_\theta(x_0,x_1) f(s_1;\theta|x_1,s_0)m(d x_1)
\\
&&\qquad \leq C  \biggl(\, \sup_{x_0 \in {\mathcal  X}} \int
p_\theta(x_0,x_1) m(d x_1) \biggr)\,d(h_1,h_2),
\end{eqnarray*}
where $0< C=\sup_{x_1 \in {\mathcal  X}} f(s_1;\theta|x_1,s_0) <
\infty$ by assumption C1 is a constant. Note that $\sup_{x_0 \in
{\mathcal  X}} \int p_\theta(x_0,x_1) m(d x_1) =1$. The equality holds
only if $h_1=h_2$ $m$-almost surely. This proves the
Lipschitz continuous condition in the second argument.

Note that C1 implies Assumption \hyperref[assK1]{K1} holds. Recall that $M_n=
{\mathbf P}(\xi_n)\circ \cdots \circ {\mathbf P}(\xi_1)\circ {\mathbf
P}(\xi_0)$ in~(\ref{5.3}). To prove the weighted mean contraction
property \hyperref[assK2]{K2}, we observe that, for $p \geq 1$,
\begin{eqnarray}\label{8.1}
\nonumber \hspace*{8mm} && \sup_{x_0,s_0} {\mathbf E}_{(x_0,s_0)}  \biggl\{ \log  \biggl( L_p
\frac{w(X_p,\xi_p)}{w(x_0,s_0)} \biggr)   \biggr\}
\\
\nonumber &&\qquad = \sup_{x_0,s_0} {\mathbf E}_{(x_0,s_0)}  \biggl\{ \log
\biggl(\, \sup_{h_1 \neq h_2} \frac{d(M_p h_1, M_p h_2)}
    {d(h_1,h_2)}\frac{w(X_p,\xi_p)}{w(x_0,s_0)} \biggr)  \biggr\}
\\
&&\qquad < \sup_{x_0,s_0} {\mathbf E}_{(x_0,s_0)}  \biggl\{ \log
\biggl( \biggl[\sup_{x_{0} \in {\mathcal  X}} \int p_\theta(x_{0},x_1)
f(\xi_1;\theta|x_1,s_0) m(d x_1)  \biggr]^p
\\
\nonumber &&\hspace*{85mm}
{}\times \frac{w(X_p,\xi_p)}{w(x_0,s_0)} \biggr) \biggr\}
\\
\nonumber &&\qquad <   0.
\end{eqnarray}
The last inequality follows from (\ref{5.1a}) in condition C1.

To verify that  Assumption \hyperref[assK3]{K3} holds,  as $m$ is
$\sigma$-finite, we have ${\mathcal  X}= \bigcup_{n=1}^{\infty} {\mathcal
X}_n$, where the ${\mathcal  X}_n$ are pairwise disjoint and $0 <
m({\mathcal  X}_n) < \infty$. Set
\begin{equation}\label{8.2}
 h(x) = \sum_{n=1}^{\infty} \frac{I_{{\mathcal  X}_n}(x)}{2^n m({\mathcal
X}_n)}.
\end{equation}
It is easy to see that $\int_{x \in {\mathcal  X}} h(x) m(dx) =1$ and,
hence, belongs to ${\mathbf M}$. Observe that
\begin{eqnarray}\label{8.3}
\nonumber && {\mathbf E} d^2\bigl({\mathbf P}(\xi_j)h,h\bigr)
\\
&&\qquad = {\mathbf E} \sup_{x_{j-1} \in {\mathcal  X}} \bigg|\int
p_\theta(x_{j-1},x_j)
\\
\nonumber &&\hspace*{30mm}
{}\times f(\xi_j;\theta|x_j,\xi_{j-1}) h(x_j) m(d x_j)-
h(x_{j-1})\bigg|.
\end{eqnarray}
By definition of $h(x)$ in (\ref{8.2}), it is piecewise constant, and
$p_\theta(x_{j-1},x_j) f(\xi_j;\break \varphi_{x_j}(\theta)|\xi_{j-1})$ is a
probability density function integrable over the subset
${\mathcal  X}_n$. These imply (\ref{8.3}) is finite.

Finally, we observe
\begin{eqnarray*}
&& \sup_{x_0,s_0} {\mathbf E}_{(x_0,s_0)}  \biggl\{ L_1 \frac{w(X_1,\xi_1)}{w(x_0,s_0)}   \biggr\}
\\
&&\qquad = \sup_{x_0,s_0} {\mathbf E}_{(x_0,s_0)}  \biggl\{ \sup_{h_1
\neq h_2} \frac{d({\mathbf P}(\xi_1)h_1,{\mathbf P}(\xi_1)h_2)}
    {d(h_1,h_2)}\frac{w(X_1,\xi_1)}{w(x_0,s_0)}  \biggr\}
\\
&&\qquad < \sup_{x_0,s_0} {\mathbf E}_{(x_0,s_0)}  \biggl\{ \sup_{x_{0}
\in {\mathcal  X}} \int p_\theta(x_{0},x_1) f(\xi_1;\theta|x_1,s_{0})
m(d x_1) \frac{w(X_1,\xi_1)}{w(x_0,s_0)}  \biggr\} <\infty.
\end{eqnarray*}
The last inequality follows from (\ref{5.1b}) in condition~C1.

Note that C5 implies the exponential moment condition of $g$. Hence,
the proof is complete.
\end{pf*}

In the proof of Lemma~\ref{L4} we omit $\theta$ for simplicity.

\begin{pf*}{Proof of Lemma~\ref{L4}}
We first prove that $\{Z_n,n \geq 0\}$ is Harris recurrent. Note that
the transition probability kernel of the Markov chain $\{(X_n,\xi_n), n
\geq 0 \}$, defined in (\ref{2.1})~and~(\ref{2.2}), has a probability
density with respect to $m \times Q$. And the iterated random functions
system, defined in (\ref{2.4})--(\ref{2.7}), also has a probability
density with respect to~$Q$. By making use the definition (\ref{3.2}),
there exists a measurable function $g\dvtx ({\mathcal  X}\times
{\mathbf R}^d \times {\mathbf M}) \times ({\mathcal  X} \times {\mathbf
R}^d \times {\mathbf M}) \to [0,\infty)$ such that
\begin{equation}\label{8.4}
{\mathbf P}(z,dz') = g(z,z')(m \times Q \times Q)(dz'),
\end{equation}
where $\int_{({\mathcal  X} \times {\mathbf R}^d) \times {\mathbf M}}
g(z,z') (m \times Q \times Q)(dz') =1$ for all $z \in ({\mathcal  X}
\times {\mathbf R}^d) \times {\mathbf M}$. For simplicity of notation,
we let $\Lambda(\cdot) :=(m \times Q \times Q) (\cdot)$ in the proof.
For given $n > 1$, let ${\mathbf P}^{n}(z,\cdot):= {\mathbf
P}_{z}(Z_{n} \in \cdot)$ for $z \in ({\mathcal  X} \times {\mathbf
R}^d) \times {\mathbf M}$. For $A \in {\mathcal  B}({\mathcal  X}
\times {\mathbf R}^d)$ and $B \in {\mathcal  B}({\mathbf M})$, define
\[
\Lambda^{n}(A \times B) :=  \int_{({\mathcal  X} \times {\mathbf R}^d)
\times {\mathbf M}}
 {\mathbf P}_{z'}\{Z_{n} \in A \times B\} \Lambda(dz').
\]
Then for all $A \in {\mathcal  B}({\mathcal  X} \times {\mathbf R}^d)$
and $B \in {\mathcal  B}({\mathbf M})$,
\begin{eqnarray*}
{\mathbf P}^{n+1}(z, A \times B) &=& \int_{({\mathcal  X} \times
{\mathbf R}^d) \times {\mathbf M}} {\mathbf P}^{n}(z', A \times B)
g(z,z') \Lambda(dz')
\\
&=& \int_{({\mathcal  X} \times {\mathbf R}^d) \times {\mathbf M}}
{\mathbf P}_{z'}\{Z_{n} \in A \times B\} g(z,z') \Lambda(dz').
\end{eqnarray*}

It is easy to see that, for given any $n > 1$, the family $({\mathbf
P}^{n+1}(z,\cdot))_{z \in ({\mathcal  X} \times {\mathbf R}^d) \times
{\mathbf M}}$ is absolutely continuous with respect to $\Lambda^{n}$.
Therefore, by the Radon--Nikodym theorem, ${\mathbf P}^n$ has a
probability density with respect to $\Lambda^n$ for all $n\geq 1$. Let
$g_{n}$ be such that
\begin{equation}\label{8.5}
 {\mathbf P}^{n+1}(z,dz') = g_{n}(z,z') \Lambda^{n}(dz'),\qquad
 z \in ({\mathcal  X}\times {\mathbf R}^d)
\times {\mathbf M},
\end{equation}
where $\int_{({\mathcal  X}\times {\mathbf R}^d) \times {\mathbf
M}}g_{n}(z,z') \Lambda^{n}(d z')=1$ for all $z \in ({\mathcal  X}\times
{\mathbf R}^d) \times {\mathbf M}$. Note that \mbox{$g_1=g$}. It is easy to
check that all $\Lambda^{n}$ are absolutely continuous with respect to~$\Pi$.

Denote $B^c$ as the complement of $B$. Since $\Pi((({\mathcal  X}\times
{\mathbf R}^d) \times {\mathbf M})^{c})=0$, also $\Lambda((({\mathcal
X}\times {\mathbf R}^d) \times {\mathbf M})^{c})=0$. Recall $g$ is
defined in (\ref{8.4}). It is obvious from the previous considerations
that we can choose $\delta>0$ sufficiently small such that
\begin{eqnarray*}
&& \int_{({\mathcal  X}\times {\mathbf R}^d) \times {\mathbf M}}
\int_{({\mathcal  X}\times {\mathbf R}^d) \times {\mathbf M}}
\int_{({\mathcal  X}\times {\mathbf R}^d) \times {\mathbf M}}
\mathbh{1}_{\{g_2 \ge \delta\}}(z_1,z_2)
\\
&&\hspace*{54mm} {}\times \mathbh{1}_{\{g \geq
\delta\}}(z_2,z_3)  \Lambda(dz_3) \Lambda^2(d z_2) \Pi(d z_1) > 0.
\end{eqnarray*}
Hence, by Lemma 4.3 of~\cite{NN86}, there exist a $\Pi$-positive set
$\Gamma_{1} \subset ({\mathcal  X}\times {\mathbf R}^d) \times {\mathbf
M}$ and a $\Lambda$-positive set $\Gamma_{2} \subset ({\mathcal
X}\times {\mathbf R}^d) \times {\mathbf M}$ such that
\[
\alpha := \inf_{z_1 \in \Gamma_{1}, z_3 \in \Gamma_{2}} \Lambda^2 \{z_2
\in ({\mathcal  X}\times {\mathbf R}^d) \times {\mathbf M}\dvtx
g_2(z_1,z_2) \ge \delta, g(z_2,z_3) \ge \delta\}  > 0.
\]
A combination of the above result with (\ref{8.4}) and (\ref{8.5})
implies
\begin{eqnarray}\label{8.6}
\nonumber {\mathbf P}^{3}(z_1, A \times B) &=&  \int_{({\mathcal
X}\times {\mathbf R}^d) \times {\mathbf M}} {\mathbf P}(z_2,A \times B)
 {\mathbf P}^2(z_1,dz_2)
\\
&\geq&  \int_{({\mathcal  X}\times {\mathbf R}^d) \times {\mathbf M}}
g_2(z_1,z_2)\int_{(A \times B) \cap \Gamma_{2}} g(z_2,z_3)  \Lambda(d
z_3) \Lambda^2(d z_2)
\\
\nonumber &\geq&  \alpha \delta^{2} \Lambda\bigl((A \times B) \cap
\Gamma_{2}\bigr)
\end{eqnarray}
for all $z_1 \in \Gamma_{1}$ and $A \times B \in {\mathcal
B}({({\mathcal  X}\times {\mathbf R}^d) \times {\mathbf M}})$.
Therefore, we obtain an absorbing set such that $\Gamma_{1}$ is a
regeneration set for $\{Z_n, n \geq 0\}$ on $({\mathcal  X}\times
{\mathbf R}^d) \times {\mathbf M}$, that is, $\Gamma_{1}$~is recurrent
and satisfies a minorization condition, namely, (\ref{8.6}). This
proves the Harris recurrence of $\{Z_n, n \geq 0\}$ on $({\mathcal
X}\times {\mathbf R}^d) \times {\mathbf M}$. Since $\{Z_n, n \geq 0\}$
possesses a stationary distribution, it is clearly positive Harris
recurrent.

Next, we give the proof of aperiodicity. If $\{Z_n, n \geq 0\}$ were
$q$-periodic with cyclic classes ${\Gamma}_{1},\ldots,{\Gamma}_{q}$,
say, then the $q$-skeleton $(Z_{nq})_{n \geq 0}$ would have stationary
distributions $\frac{\Pi(\cdot \cap {\Gamma}_{k})}{\Pi({\Gamma}_{k})}$
for $k=1,\ldots,q$. On the other hand, $Z_{qn}$ is aperiodic by
definition, and $M_{nq}$ is also a Markovian iterated random functions
system of Lipschitz maps, satisfying condition C1, and thus possesses
only one stationary distribution. Consequently, $q=1$ and $\{Z_n, n
\geq 0\}$ is aperiodic. Since the Markov chain $\{((X_n,\xi_n),M_n), n
\geq 0\}$ has a probability density with respect to $\Lambda$, it is
obviously $\Lambda$-irreducible. The proof is complete.
\end{pf*}

\begin{pf*}{Proof of Lemma~\ref{L5}}
In order to define the Fisher information (\ref{5.6}), we need to
verify that there exists a $\delta > 0,$ such that ${\partial  \log \|
{\mathbf P}_{\theta}(\xi_1)\circ {\mathbf
P}_{\theta}(\xi_0)\pi\|}/{\partial \theta} \in L_2({\mathbf
P}_{\Pi}^{\theta})$ for $\theta \in N_\delta(\theta_0),$ a
$\delta$-neighborhood of $\theta_0$. That is, we need to show
\begin{equation}\label{8.7}
 {\mathbf E}^{\theta}_{\Pi}  \biggl(\frac{\partial  \log
\| {\mathbf P}_{\theta}(\xi_1)\circ {\mathbf
P}_{\theta}(\xi_0)\pi\|}{\partial \theta}  \biggr)^2 < \infty,
\end{equation}
for $\theta \in N_\delta(\theta_0).$

It is easy to see that C5 implies that
\[
\sup_{x \in {\mathcal  X}} E^{\theta}_x  \biggl(\frac{\partial \log
\int_{y \in {\mathcal  X}} \pi(x) p(x,y)
f(\xi_0;\theta|x)f(\xi_1;\theta|y,\xi_0)m(dy)} {\partial \theta}
\biggr)^2 < \infty
\]
for $\theta \in N_\delta(\theta_0).$ And this leads to
\begin{equation}\label{8.8}
\hspace*{10mm}  \sup_{x \in {\mathcal  X} } {\mathbf E}^{\theta}_{x}
\biggl(\frac{\partial  \log \int_{y \in {\mathcal  X}} \pi(x)
p(x,y)f(\xi_0; \theta|x) f(\xi_1;\theta|y,\xi_0)m(dy)}{\partial \theta}
\biggr)^2 < \infty
\end{equation}
for $\theta \in  N_\delta(\theta_0),$ where ${\mathbf E}_x^{\theta}$ is
the expectation under ${\mathbf P}^{\theta} (\cdot,\cdot)$.

Finally, (\ref{8.8}) implies (\ref{8.7}) and we have the proof.
\end{pf*}

\begin{pf*}{Proof of Lemma \ref{L6}}
For each $j=1,\ldots, q$,
\begin{eqnarray*}
\frac1{\sqrt{n}}l'_j(\theta_0)&=& \frac1{\sqrt{n}}
\frac{\partial}{\partial\theta_j} \log \| {\mathbf P}_\theta(\xi_n)
\circ \cdots \circ {\mathbf P}_\theta(\xi_1) \circ
{\mathbf P}_\theta(\xi_0) \pi \|\bigg|_{\theta=\theta_0}
\\
&=& \frac1{\sqrt{n}}\sum_{k=1}^{n} \biggl(
\frac{\partial}{\partial\theta_j} \log  \frac{\| {\mathbf
P}_\theta(\xi_k) \circ \cdots \circ {\mathbf P}_\theta(\xi_1) \circ
{\mathbf P}_\theta(\xi_0) \pi \|}{\| {\mathbf P}_\theta(\xi_{k-1})
\circ \cdots \circ {\mathbf P}_\theta(\xi_1) \circ {\mathbf
P}_\theta(\xi_0) \pi \|}\bigg|_{\theta=\theta_0}  \biggr)
\\
&=& \frac{1}{\sqrt{n}}\sum_{k=1}^{n}\frac{\partial}{\partial\theta_j}
 g(M_{k-1},M_k)\bigg|_{\theta= \theta_0}.
\end{eqnarray*}

Now, for each $h \in {\mathbf M}$, $\alpha = (\alpha_1,\ldots,\alpha_q)
\in C^q$, and a $({\mathcal  X} \times {\mathbf R}^d) \times {\mathbf
M}$ measurable function $\varphi$ with $\|\varphi\|_{wh} < \infty$,
define
\begin{eqnarray}\label{8.9}
\nonumber \hspace*{4mm} &&  ({\mathbf T}_1(\alpha)\varphi)\bigl((x,s), h\bigr)
\\
&&\qquad = {\mathbf E}^{\theta_0}_{(x,s)} \biggl\{ \exp  \biggl(
(\alpha_1,\ldots,\alpha_q)' \biggl(
\frac{\partial}{\partial\theta_1}\log \|{\mathbf P}_\theta(\xi_1) \circ
{\mathbf P}_\theta(\xi_0) h\|\bigg|_{\theta=\theta_0},\ldots,
\nonumber
\\[-8pt]
\\[-8pt]
\nonumber
&&\hspace*{6.3mm} \phantom{\qquad = {\mathbf E}^{\theta_0}_{(x,s)} \biggl\{ \exp
\biggl( (\alpha_1,\ldots,\alpha_q)' \biggl(}
{}\frac{\partial}{\partial\theta_q}\log \|{\mathbf P}_\theta(\xi_1)
\circ {\mathbf P}_\theta(\xi_0) h\|\bigg|_{\theta=\theta_0}\biggr)
\biggr)
\\
\nonumber &&\hspace*{59mm}
{}\times \varphi\bigl((X_1,\xi_1),{\mathbf P}_\theta(\xi_1) \circ
{\mathbf P}_\theta(\xi_0) h(x)\bigr)  \biggr\}.
\end{eqnarray}
By using an argument similar to that of Lemma \ref{L2}, we have, for
sufficiently small~$|\alpha|$, ${\mathbf T}_1(\alpha)$ is a bounded and
analytic operator. Let $\lambda_{{\mathbf T}_1}^{\theta_0}(\alpha)$ be
the eigenvalue of ${\mathbf T}_1(\alpha)$ corresponding to a
one-dimensional eigenspace. Define $\gamma_j$ as that in
Lemma~\ref{L2}(v). By conditions C1--C5 and Lemma \ref{L4}, it is easy to see
that
\begin{equation}\label{8.10}
\hspace*{8mm} \gamma_j = \frac{\partial}{\partial \alpha_j} \lambda_{{\mathbf
T}_1}^{\theta_0} (\alpha)\bigg|_{\alpha =0} =  {\mathbf
E}^{\theta_0}_\Pi \biggl(\frac{\partial}{\partial \theta_j} \log \|
{\mathbf P}_\theta(\xi_1) \circ {\mathbf P}_\theta(\xi_0)
\pi\|\bigg|_{\theta =\theta_0} \biggr)= 0.
\end{equation}
By Corollary~\ref{C1}, we have
\begin{equation}\label{8.11}
 \frac1{\sqrt{n}} (l'_j(\theta_0) )_{j=1,\ldots,q}
 \longrightarrow N\bigl(0,{\bolds\Sigma}(\theta_0)\bigr)\qquad\mbox{in distribution,}
\end{equation}
where the variance--covariance matrix
\begin{equation}\label{8.12}
{\bolds\Sigma}(\theta_0) = (\Sigma_{ij}(\theta_0))=
\biggl(\frac{\partial^2 \lambda_{{\mathbf T}_1}^{\theta_0}(\alpha)}{
\partial \alpha_i \,\partial \alpha_j}\bigg|_{\alpha=0 }
 \biggr)_{i,j=1,\ldots,q}.
\end{equation}

In the following, we will verify that the variance--covariance matrix
${\bolds\Sigma}(\theta_0)$ defined as (\ref{8.12}) is the Fisher
information matrix ${\mathbf I}(\theta_0)$. By Lemma \ref{L2} and
Corollary~\ref{C1}, we have
\begin{eqnarray*}
&&  {\mathbf E}^{\theta_0}_\Pi  \biggl(\biggl(\frac{\partial}{\partial\theta_j}\log\|M_n
 \pi \| \bigg|_{\theta=\theta_0}\biggr)\biggl(\frac{\partial}{\partial\theta_k}\log\|M_n \pi
 \| \bigg|_{\theta=\theta_0}\biggr) \biggr) -  n \frac{\partial^2}{\partial \alpha_j\, \partial
 \alpha_k} \lambda_{{\mathbf T}_1}^{\theta_0}(\alpha)\bigg|_{\alpha = 0}
\\
&&\qquad  \longrightarrow 0
\end{eqnarray*}
as $n \rightarrow \infty$. Therefore,
\begin{eqnarray*}
  \Sigma_{jk}(\theta_0)
&=&  \frac{\partial^2}{\partial \alpha_j \,\partial
 \alpha_k} \lambda_{{\mathbf T}_1}^{\theta_0}(\alpha)\bigg|_{\alpha = 0}
\\
&=&  \lim_{n\rightarrow\infty}\frac{1}{n}
 {\mathbf E}^{\theta_0}_\Pi  \biggl(\frac{\partial}
 {\partial\theta_j}\log \|M_n \pi \| \bigg|_{\theta=\theta_0} \biggr)
 \biggl(\frac{\partial} {\partial\theta_k}\log \|M_n \pi\| \bigg|_{\theta=\theta_0} \biggr)
\\
&=& \lim_{n \rightarrow \infty} - \frac{1}{n} {\mathbf
E}_\Pi^{\theta_0}
   \biggl(\frac{\partial^2}{\partial\theta_j\, \partial \theta_k}
 \log \|M_n \pi \| \bigg|_{\theta=\theta_0} \biggr)
\\
&=&   -  {\mathbf E}_\Pi^{\theta_0}
   \biggl(\frac{\partial^2}{\partial\theta_j\, \partial \theta_k}
 \log \| {\mathbf P}_\theta(\xi_1) \circ {\mathbf P}_\theta(\xi_0)
\pi \| \bigg|_{\theta=\theta_0} \biggr)
\\
&=& {\mathbf E}_\Pi^{\theta_0}  \biggl(
\frac{\partial}{\partial\theta_j}\log\| {\mathbf P}_\theta(\xi_1) \circ
{\mathbf P}_\theta(\xi_0) \pi \| \bigg|_{\theta=\theta_0} \biggr)
\\
&&{}\times \biggl( \frac{\partial}{\partial\theta_k} \log\|{\mathbf
P}_\theta(\xi_1) \circ {\mathbf P}_\theta(\xi_0) \pi \|
\bigg|_{\theta=\theta_0} \biggr)
\\
&=&  I_{jk}(\theta_0).
\end{eqnarray*}\upqed
\end{pf*}

\begin{appendix}
\section*{Appendix}\label{app}
\renewcommand{\thelemma}{A.\arabic{lemma}}
\setcounter{lemma}{0}
\setcounter{equation}{0}

\subsection*{\texorpdfstring{Proofs of Lemma \protect\ref{L1} and
Theorem~\protect\ref{T2}}{Proofs of Lemma 1 and Theorem~2}}
In the following proofs we will use the same notation as in
Sections \ref{s3} and \ref{s4} unless specified. Without loss of
generality, in this section we consider the case $M_0=\mbox{Id}$, the
identity, and the transition probability ${\mathbf P}$ of the Markov
chain $\{(Y_n,M_n), n \geq 0\}$ depends on the initial state $Y_0=y$
only. Denote it as ${\mathbf P}_y$, and let ${\mathbf E}_y$ be the
corresponding expectation. To prove Lemma \ref{L1}, we need the following
lemma first.

\begin{lemma}\label{L7}
Let $\{(Y_n,M_n),n\geq 0\}$ be the MIRFS of Lipschitz functions defined
in~\textup{(\ref{2.1})} satisfying Assumption~\textup{\ref{assK}}.
There exists $0< \delta_0 < 1$ such that, for all $0< \delta \leq
\delta_0$, there exist $K > 0$, and $0 < \eta< 1$, so that
\[
\sup_y {\mathbf E}_y  \biggl\{  \biggl( \frac{d(M_n^u, M_n^v)}{d(u, v)}
\frac{w(Y_n)}{w(y)} \biggr)^{\delta} \biggr\} \leq K
\eta^n,\qquad\mbox{for }n \in N\mbox{ and }u,v \in {\mathbf M}.
\]
\end{lemma}

\begin{pf}
For given $0< \delta < 1$, and $y \in {\mathcal  Y}$, denote
\[
c_n(y)=\sup  \biggl\{ {\mathbf E}_y  \biggl[  \biggl(\frac{d(M_n^u,
M_n^v)}{d(u,v)} \frac{w(Y_n)}{w(y)} \biggr)^{\delta} \biggr]\dvtx
 u, v \in \mathbf{M}  \biggr\},
\]
and let $\eta_n= \sup\{c_n(y),y\in {\mathcal  Y}\}.$ Denote $u_m
=M_m^u$ and $v_m=M_m^v$. Let ${\mathcal  F}_m$ be the $\sigma$-algebra
generated by $\{(Y_k,M_k), 0 \leq k \leq m\}$. Then
\begin{eqnarray*}
&&  {\mathbf E}_y \biggl\{  \biggl(\frac{d(M_{n+m}^u,M_{n+m}^v)}{d(u,
v)} \frac{w(Y_{n+m})}{w(y)} \biggr)^{\delta}\Big|{\mathcal  F}_m
\biggr\}
\\
&&\qquad = {\mathbf E}_y  \biggl\{  \biggl(\frac{d(F_{n:m}(M_m^u),
F_{n:m}(M_m^v))}{d(u, v)}\frac{w(Y_{n+m})}{w(y)} \biggr)^{\delta} \Big|{\mathcal  F}_m \biggr\}
\\
&&\qquad =  \biggl( \frac{d(M_m^u,M_m^v)}{d(u,v)} \frac{w(Y_m)}{w(y)}
\biggr)^{\delta}    {\mathbf E}_y \biggl\{
\biggl(\frac{d(F_{n:m}(u_m),F_{n:m}(v_m))}{d(u_m,v_m)} \frac{w(Y_{n+m})}{w(Y_m)} \biggr)
^{\delta} \Big|{\mathcal  F}_m \biggr\}
\\
&&\qquad =  \biggl( \frac{d(M_m^u,M_m^v)}{d(u,v)} \frac{w(Y_m)}{w(y)}
\biggr)^{\delta}     {\mathbf E}_{Y_m}  \biggl\{
\biggl(\frac{d(M_{n}^{u_m},M_{n}^{v_m})}{d(u_m,v_m)}
\frac{w(Y_{n+m})}{w(Y_m)} \biggr)^{\delta}  \biggr\}
\\
&&\qquad \leq  \biggl(\frac{d(M_m^u,M_m^v)}{d(u,v)} \frac{w(Y_m)}{w(y)}
\biggr)^{\delta} c_n(Y_m) \leq \eta_n
\biggl(\frac{d(M_m^u,M_m^v)}{d(u,v)}\frac{w(Y_m)}{w(y)}
\biggr)^{\delta}.
\end{eqnarray*}
This implies that
\begin{eqnarray*}
&& {\mathbf E}_y  \biggl\{  \biggl(\frac{d(M_{n+m}^u,M_{n+m}^v)}{d( u,
v)}\frac{w(Y_{n+m})}{w(y)} \biggr)^{\delta}  \biggr\}
\\
&&\qquad \leq \eta_n
{\mathbf E}_y \biggl\{
\biggl(\frac{d(M_m^u,M_m^v)}{d(u,v)}\frac{w(Y_m)}{w(y)}
\biggr)^{\delta} \biggr\},
\end{eqnarray*}
or $\eta_{n+m} \leq \eta_n \eta_m.$ Therefore,
\begin{equation}\label{7.1}
\lim_{n \to \infty}\eta^{1/n}_n= \inf \{\eta^{1/n}_n,n\in N\}.
\end{equation}

It is known by Assumption \hyperref[assK2]{K2} that there exist $p \geq
1$ and $d > 0$ such that $ \sup_y {\mathbf E}_y  \{ \log  (
\frac{d(M_p^u, M_p^v)}{d(u, v)}\frac{w(Y_{p})}{w(y)} )  \} < -d < 0.$
Along with $\sup_y {\mathbf E}_y \{ \frac{w(Y_p)}{w(y)}\}< \infty$
by~(\ref{4.2}) and $\sup_y {\mathbf E}_y \{ \frac{l(F_1)w(Y_1)}{w(y)}\}<
\infty$ by Assumption \hyperref[assK3]{K3}, we have
\[
\eta_p \leq \sup_{y \in {\mathcal  Y}} {\mathbf E}_y  \biggl\{
\biggl(l(F_1)^p \frac{w(Y_p)}{w(y)} \biggr)^{\delta}  \biggr\}
:=
\sup_{y \in {\mathcal Y}} {\mathbf E}_y  \biggl\{ \exp  \biggl( \delta
G_p + \delta \log \frac{w(Y_p)}{w(y)} \biggr) \biggr\} < \infty,
\]
where $G_p = p \log l(F_{1})$.

Since $e^y \leq 1+y+y^2 e^{|y|}/2,$ we have, for $y \in {\mathcal  Y},
u,v \in {\mathbf M},$
\begin{eqnarray*}
&& {\mathbf E}_y \biggl\{
\biggl(\frac{d(M_p^u,M_p^v)}{d(u,v)}\frac{w(Y_{p})}{w(y)}
\biggr)^{\delta} \biggr\}
\\
&&\qquad \leq 1  + \delta {\mathbf E}_y  \biggl\{\log
\biggl(\frac{d(M_p^u,M_p^v)}{d(u, v)}\frac{w(Y_p)}{w(y)} \biggr) \biggr\}
\\
&&\qquad\quad{}  +  \delta^2 {\mathbf E}_y  \biggl\{  \biggl(G_p + \log
\frac{w(Y_p)}{w(y)} \biggr)^2 \exp \biggl( \delta G_p +  \delta \log
\frac{w(Y_p)}{w(y)} \biggr)  \biggr\}.
\end{eqnarray*}
For $u,v \in {\mathbf M}$, we have
\[
\eta_p \leq 1-d \delta +  \delta^2 \sup_{y \in {\mathcal  Y}} {\mathbf
E}_y  \biggl\{  \biggl(G_p + \log \frac{w(Y_p)}{w(y)} \biggr)^2 \exp
\biggl( \delta G_p + \delta \log \frac{w(Y_p)}{w(y)} \biggr)  \biggr\}.
 \]
Therefore, we can choose $\delta_0 > 0$ small enough so that $\eta_p <
1$. Along  with (\ref{7.1}), we obtain the proof.
\end{pf}

\begin{pf*}{Proof of Lemma~\ref{L1}}
For given $\varphi \in {\mathcal  H}$, $y \in {\mathcal  Y},$ and $u,v
\in {\mathbf M},$ if $m \leq n$, we have, for $0 < \delta \leq \delta_0
< 1$,
\begin{eqnarray*}
&&  \big|{\mathbf T}^n\varphi(y,u) - {\mathbf E}_y \varphi\bigl(Y_n, F_{n:m}(v)\bigr)\big|/w(y)
\\
&&\qquad = \big|{\mathbf E}_y\varphi(Y_n,M_n^u) - {\mathbf E}_y \varphi\bigl(Y_n, F_{n:m}(v)\bigr)\big|/w(y)
\\
&&\qquad \leq \|\varphi\|_h {\mathbf E}_y \bigl\{ d\bigl(M_n^u, F_{n:m}(v)\bigr)^\delta w(Y_n)^\delta \bigr\}/w(y)
\\
&&\qquad \leq \|\varphi\|_h {\mathbf E}_y  \biggl\{ {\mathbf E}_y
\biggl[ \biggl(d\bigl(F_{n:m}(M_{n-m}^u), F_{n:m}(v)\bigr)
\frac{w(Y_n)}{w(Y_{n-m})} \biggr)^\delta \Big|
{\mathcal  F}_{n-m}  \biggr]  \frac{w(Y_{n-m})^\delta}{w(y)}  \biggr\}
\\
&&\qquad \leq \|\varphi\|_h {\mathbf E}_y  \biggl\{\sup_{u,v \in
{\mathbf M}} {\mathbf E}_{Y_{n-m}} \biggl[  \biggl(d(M_{m}^u, M_{m}^v)
\frac{w(Y_n)}{w(Y_{n-m})} \biggr)^\delta \biggr]
\frac{w(Y_{n-m})^\delta}{w(y)} \biggr\}
\\
&&\qquad \leq \|\varphi\|_h {\mathbf E}_y  \biggl\{\sup_{u,v \in
{\mathbf M}} {\mathbf E}_{Y_{n-m}} \biggl[  \biggl(\frac{d(M_{m}^u,
M_{m}^v)}{d(u,v)} \frac{w(Y_n)}{w(Y_{n-m})} \biggr)^\delta  \biggr]
\frac{w(Y_{n-m})}{w(y)}    \biggr\}.
\end{eqnarray*}
Note that in the last inequality we use $d(u,v) \leq 1$ and $w(y) \geq
1$ for all $y \in {\mathcal  Y}$.

By making use of Lemma \ref{L7}, and $\sup_{y \in {\mathcal  Y}} E_y
[w(Y_1)/w(y)] < \infty$ in (\ref{4.2}), there exist  $K > 0$ and $ 0<
\eta < 1$ such that
\begin{equation}\label{7.2}
\hspace*{12pt}
\big|{\mathbf T}^n\varphi(y,u) - {\mathbf E}_y \varphi\bigl(Y_n,
F_{n:m}(v)\bigr)\big|/w(y) \leq \|\varphi\|_h K \eta^m \leq
\|\varphi\|_{wh} K \eta^m.
\end{equation}

Denote $h(y)= {\mathbf E}_y \varphi(Y_m, F_{m}(v))$. Then by assumption
(\ref{4.1}), there exist $\gamma > 0$ and $0 < \rho <1$ such that
\begin{eqnarray}\label{7.3}
\nonumber && \big|{\mathbf E}_y \varphi\bigl(Y_n, F_{n:m}(v)\bigr)-
{\mathbf E}_{\Pi} \varphi\bigl(Y_m, F_{m}(v)\bigr)\big|/w(y)
\\
&&\qquad \leq \big|{\mathbf E}_y\bigl\{ {\mathbf E}_{Y_{n-m}}
\varphi\bigl(Y_m, F_m(v)\bigr)\bigr\}
- {\mathbf E}_{\Pi} \varphi\bigl(Y_m, F_{m}(v)\bigr)\big|/w(y)
\nonumber
\\[-8pt]
\\[-8pt]
\nonumber
&&\qquad \leq \bigg|{\mathbf E}_y h(Y_{n-m}) - \int h(y)
\Pi(dy)\bigg|\Big/w(y)
\\
\nonumber &&\qquad \leq \|\varphi\|_{wh} \gamma \rho^{n-m}.
\end{eqnarray}

For given $m,k \in N$, by using Lemma \ref{L7} again we have
\begin{eqnarray*}
&&  \big|{\mathbf E}_{\Pi} \varphi\bigl(Y_m, F_{m}(v)\bigr)
- {\mathbf E}_{\Pi} \varphi\bigl(Y_{m+k}, F_{m+k}(v)\bigr)\big|/w(y)
\\
&&\qquad \leq {\mathbf E}_{\Pi}  \bigl\{ \big|
\varphi\bigl(Y_{m+k}, F_{m+k:k}(v)\bigr) -
\varphi\bigl(Y_{m+k}, F_{m+k:k}(M_k^v)\bigr)\big|\bigr\}/w(y)
\\
&&\qquad \leq \|\varphi\|_h {\mathbf E}_{\Pi}  \biggl\{ d\bigl(F_{m+k:k}(v),
F_{m+k:k}(M_k^v)\bigr)^{\delta} \frac{w(Y_{m+k})^\delta}{w(y)}  \biggr\}
\\
&&\qquad \leq \|\varphi\|_h {\mathbf E}_{\Pi} \biggl\{ \sup_{u,v \in {\mathbf M}}
{\mathbf E}_{Y_m} \biggl[  \biggl(\frac{d(M_{m}^u, M_{m}^v)}{d(u,v)}
\frac{w(Y_{m+k})}{w(Y_m)} \biggr)^{\delta} \biggr] \frac{w(Y_m)}{w(y)}  \biggr\}
\\
&&\qquad \leq \|\varphi\|_{wh} K \eta^m.
\end{eqnarray*}

By making use of (\ref{7.2}), (\ref{7.3}) and the above inequality, we
have that for any given $n \geq m$, $k \geq 0$, and for all $u,v \in
{\mathbf M}$,
\[
\big|{\mathbf T}^n\varphi(y,u) - {\mathbf E}_{\Pi} \varphi\bigl(Y_{m+k}, F_{m+k}(v)\bigr)\big|/w(y)
\leq \|\varphi\|_{wh}(2 K \eta^m + \gamma \rho^{n-m}).
\]
By setting $m = n/2$, we have that there exist $A >0$ and $0< r < 1$
such that
\begin{equation}\label{7.4}
 \|{\mathbf T}^n\varphi(y,u) - {\mathbf Q}\varphi(y,u)\|_w
\leq \|\varphi\|_{wh} Ar^{n}.
\end{equation}

On the other hand, for $u,v \in {\mathbf M}$,
\begin{eqnarray}\label{7.5}
\nonumber && \frac{|({\mathbf T}^n-{\mathbf Q})\varphi(y,u) -
({\mathbf T}^n-{\mathbf Q})\varphi(y,v)|}{(w(y)\,d(u,v))^{\delta}}
\\
\nonumber &&\qquad =   \bigg|{\mathbf E}_y\varphi(Y_n,M_n^u) - \int \varphi(y,u)
    \Pi(dy \times du)
\\
\nonumber &&\hspace*{12mm} {}
- {\mathbf E}_{y}\varphi(Y_n,M_n^v) + \int \varphi(y,v)
    \Pi(dy \times dv) \biggr|
\\
&&\qquad\quad{} \times
\bigl[{\bigl(w(y)\,d(u,v)\bigr)^{-\delta}}\bigr]^{-1}
\\
\nonumber &&\qquad \leq  \frac{ {\mathbf E}_y\{  |\varphi(Y_n,M_n^u)- \varphi(Y_n,M_n^v) |\}}{(w(y)\,d(u,v))^{\delta}}
\\
\nonumber &&\qquad \leq \|\varphi\|_h  \sup_{y}{\mathbf E}_y \biggl\{
\biggl(\frac{d(M_n^u,
M_n^v)}{d(u,v)}\frac{w(Y_{n})} {w(y)} \biggr)^{\delta} \biggr\}
\\
\nonumber &&\qquad \leq
\|\varphi\|_{wh} K \eta^n\qquad \mbox{by Lemma~\ref{L7}}.
\end{eqnarray}

Denote $\rho_*=\min\{\eta, r\}$ and $\gamma_*= A+K$. Combine
(\ref{7.4}) and (\ref{7.5}) to get
\[
\|{\mathbf T}^n  - {\mathbf Q}  \|_{wh} = \sup_{\varphi \in {\mathcal
H}, \|\varphi \|_{wh} \leq 1} \|{\mathbf T}^n \varphi - {\mathbf Q}
\varphi \|_{wh}
 \leq \sup_{\varphi \in {\mathcal  H}, \|\varphi \|_{wh} \leq 1}
\|\varphi\|_{wh} \gamma_* \rho_*^n \leq \gamma_*  \rho_*^n.
\]
Then we have (\ref{4.15}) and this completes the proof.
\end{pf*}

\begin{pf*}{Proof of Theorem~\ref{T2}}
By using Lemma \ref{L2}, standard arguments
involving smoothing inequalities and Fourier inversion (cf. Chapter~4
of~\cite{BR76}) reduce the proof to that of
showing for every $\delta > 0, a > 0$ and $b > 1,$
\begin{equation}\label{7.11}
\sup_{\delta\leq|\alpha|\leq n^a}\big|E_{\pi} \bigl(e^{i\alpha' S_n} \bigr) \big| =
o(n^{-b}).
\end{equation}

To prove (\ref{7.11}), we follow the same idea as (3.43) of~\cite{GH83},
letting $\zeta_t =S_t-S_{t-1}\ (t=1,2,\ldots),
\zeta_0= S_0$ and $\tilde{\varphi}((y,u),(y',v)) =
E\{e^{i\alpha'\zeta_1}|(Y_0=y,M_0=u),(Y_1=y',M_1=v)\}$.

Let $J = \{1,\ldots,n\}$, and fix $m > 1$ to be determined later.
Divide $J$ into blocks $A_{1},B_{1},\ldots,A_{l},B_{l}$ as follows.
Define $j_{1},\ldots,j_{l}$ by $j_{1}= 1,$ and $j_{k+1}=\inf \{j \geq
j_{k}+7m\dvtx j\in J \},$ and let $l$ be the smallest integer for which the
$\inf$ is undefined. Write
\begin{eqnarray*}
A_{k} &=&   \prod \bigl\{e^{n^{-1/2} i \alpha'  \zeta_{j}}\dvtx |j-j_{k}| \leq m\bigr\},\hspace*{30mm}\qquad k=1,\ldots,l,
\\
B_{k} &=&   \prod \bigl\{e^{n^{-1/2} i \alpha'  \zeta_{j}}\dvtx j_{k}+m+1 \leq j
\leq j_{k+1}-m-1 \bigr\},\qquad  k=1,\ldots,l-1,
\\
B_{l}&=&  \prod \bigl\{e^{n^{-1/2} i \alpha'  \zeta_{j}}\dvtx j > j_{l}+m+1 \bigr\}.
\end{eqnarray*}
Then $e^{i\alpha' S_{n}}=   \prod_{k=1}^{l} A_{k}B_{k}.$ Given $y \in
{\mathcal  Y}$, we have
\begin{eqnarray}\label{7.12}
&&  \Bigg|E_y   \prod_{1}^{l} A_{k}B_{k} - E_y   \prod_{1}^{l}B_{k}E(A_{k}|\zeta_{j}\dvtx j\neq j_{k}) \Bigg|
\nonumber
\\[-8pt]
\\[-8pt]
\nonumber
&&\qquad \leq \sum_{q=1}^{l}  \Bigg| E_y    \prod_{1}^{q-1} A_{k}B_{k}
\bigl(A_{q}-E(A_{k}|\zeta_{j}\dvtx j\neq j_{q})\bigr)
\prod_{q+1}^{l}B_{k}E(A_{k}|\zeta_{j}\dvtx j\neq j_{k})
\Bigg|.\hspace*{-8mm}
\end{eqnarray}
By using Lemma \ref{L2}(iv), there exists $\delta > 0$ such that $ E|
E(A_{k}|\zeta_{j}\dvtx j\neq j_{q})- E(A_{k}|\zeta_{j}\dvtx 0<|j-j_{k}| \leq  3m
)| \leq e^{-\delta m}.$ Therefore, (\ref{7.12}) $\leq$
\begin{eqnarray}\label{7.13}
\nonumber &&  \sum_{q=1}^{l}  \Bigg| E_y    \prod_{1}^{q-1} A_{k}B_{k} \bigl(A_{q}-E(A_{k}|\zeta_{j}\dvtx
j\neq j_{q})\bigr)
\\
&&\hspace*{6mm}
{}\times   \prod_{q+1}^{l}B_{k}E(A_{k}|\zeta_{j}\dvtx 0<|j-j_{k}| \leq  3m )  \Bigg|
\\
\nonumber
&&\qquad {}   +\sum_{q=1}^{l} e^{-\delta m}.
\end{eqnarray}

The first summation term in (\ref{7.13}) vanishes since $
\prod_{1}^{q-1}A_{k}B_{k}$ and $\prod_{q+1}^{l}B_{k}\times E(A_{k}|\zeta_{j}\dvtx 0<|j-j_{k}| \leq 3m)$ are both
measurable with respect to the $\sigma$-field generated by
$\zeta_{j}:j\neq j_{q}$.

Recall that the functions $E(A_{k}|\zeta_{j}\dvtx 0<|j-j_{k}| \leq 3m),
\mbox{for } k=1,\ldots,l$, are weakly dependent since $j_{k+1}-j_{k}\geq
7m, k=1,\ldots,l-1$. Using Assumption \hyperref[assK1]{K1},
(\ref{4.18}) and (\ref{4.19}), we obtain
\begin{eqnarray*}
&& \Bigg| E_y   \prod_{1}^{l}B_{k}E(A_{k}|\zeta_{j}\dvtx 0<|j-j_{k}|\leq 3m) \Bigg|
\\
&&\qquad \leq E_y  \Bigg|   \prod_{1}^{l}E(A_{k}| \zeta_{j}\dvtx 0<|j-j_{k}| \leq 3m)  \Bigg|
\\
&&\qquad \leq    \prod_{1}^{l} E_y \big|E(A_{k}| \zeta_{j}\dvtx 0<|j-j_{k}| \leq 3m)\big|+l
e^{-\delta m}.
\end{eqnarray*}
With the strong nonlattice condition (\ref{4.20}), and conditional
strong nonlattice condition (\ref{4.21}), we find an upper bound for
$E_y|E(A_{k}|\zeta_{j}\dvtx 0<|j-j_{k}|\leq 3m)|.$

We have for $|\alpha| \geq \delta$ the relation $E_y|E(A_{k}|
\zeta_{j}\dvtx j\neq j_{q})| \leq e^{-\delta}$ and, hence, by~(\ref{4.21})
for all $\alpha \in {\mathbf R}^p,  | \alpha | \leq \delta$, $
E_y|E(A_{k}|\zeta_{j}\dvtx j\neq j_{q})| \leq \exp(- \delta | \alpha
|^{2}/n).$ Therefore, for all $\alpha \in {\mathbf R}^p$,
\begin{eqnarray*}
&&   E_y\big|E(A_{k}|\zeta_{j}\dvtx 0<|j-j_{k}|\leq 3m)\big|
\\
&&\qquad \leq e^{-\delta m}+ E_y\big|E(A_{k}| \zeta_{j}\dvtx  j\neq j_{q})\big| \leq
e^{-\delta m} + \max \bigl(\exp(-\delta |\alpha|^{2}/n),e^{-\delta}\bigr).
\end{eqnarray*}

If we choose $K$ appropriately and let $m$ be the integral part of $K
\log n$, then the assertion of the lemma follows from $\exp (- \delta
|\alpha|^{2}/n)^{{n}/{m}} \leq \exp(-\delta |\alpha|^{2}/\break (K \log n))
\leq \exp(-\delta^{'}n^{\varepsilon/2})$ for $|\alpha| \geq
cn^{\varepsilon}$ and some $\delta^{'}>0$.
\end{pf*}
\end{appendix}

\section*{Acknowledgments}
The author is grateful to the Editor Professor Jianqing Fan,
an Associate Editor and a referee for constructive comments, suggestions
and correction of some errors in the earlier version.

\printaddresses

\end{document}